\documentclass[11pt]{article}%
\usepackage{amsfonts}
\usepackage{amsmath}
\usepackage{amssymb}
\usepackage{amsxtra}
\usepackage{graphicx}
\usepackage{geometry}
\usepackage{color}
\usepackage{colortbl}%
\providecommand{\U}[1]{\protect\rule{.1in}{.1in}}
\newtheorem{theorem}{Theorem}

\newtheorem{remark}[theorem]{Remark}

\numberwithin{equation}{section}
\ifx\pdfoutput\relax\let\pdfoutput=\undefined\fi
\newcount\msipdfoutput
\ifx\pdfoutput\undefined\else
\ifcase\pdfoutput\else
\ifx\paperwidth\undefined\else
\ifdim\paperheight=0pt\relax\else\pdfpageheight\paperheight\fi
\ifdim\paperwidth=0pt\relax\else\pdfpagewidth\paperwidth\fi
\fi\fi\fi

\begin{document}

\title{Spectra of open waveguides in periodic media.}
\author{G.Cardone\thanks{Universit\`{a} del Sannio, Department of Engineering, Corso
Garibaldi, 107, 82100 Benevento, Italy; email: giuseppe.cardone@unisannio.it.}%
, S.A.Nazarov\thanks{Mathematics and Mechanics Faculty, St. Petersburg State
University, 198504, Universitetsky pr., 28, Stary Peterhof, Russia;
Saint-Petersburg State Polytechnical University, Polytechnicheskaya ul., 29,
St. Petersburg, 195251, Russia; Institute of Problems of Mechanical
Engineering RAS, V.O., Bolshoj pr., 61, St. Petersburg, 199178, Russia; email:
srgnazarov@yahoo.co.uk.}, J. Taskinen\thanks{University of Helsinki, Department of Mathematics and Statistics,
P.O. Box 68, 00014 Helsinki, Finland;
email: jari.taskinen@helsinki.fi.}}
\maketitle

\begin{abstract}
We study the essential spectra of formally self-adjoint elliptic systems on
doubly periodic planar domains perturbed by a semi-infinite periodic row of foreign inclusions. We show that
the essential spectrum of the problem consists of the essential spectrum of
the purely periodic problem and another component, which is the union of the
discrete spectra of model problems in the infinite perturbation strip; these
model problems arise by an application of the partial Floquet-Bloch-Gelfand transform.

\medskip

Keywords: Spectral bands, double periodic medium, unbounded periodic
perturbation, open waveguide

\medskip

MSC: 35P05, 47A75, 49R50, 78A50

\end{abstract}

\section{Introduction.\label{sect1}}

\subsection{Preamble.\label{sect1.1}}

\begin{figure}[ptb]
\begin{center}
\includegraphics[height=1.5in,
width=4in
]
{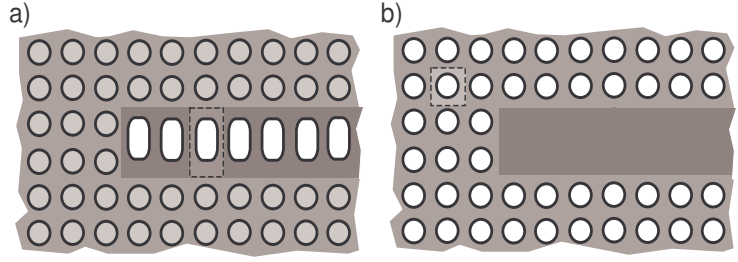}
\end{center}
\caption{a) Domain $\Omega$ (grey) including the open $\mathsf{I}$-shaped  waveguide
$\Omega^{1+}$ (dark grey). b) Strip $\varpi^+$ (dark grey), which
coincides with the open waveguide
$\Omega^{1+}$ in the special case $\vartheta^{1}=\varnothing$.
\hfill \break Dotted lines indicate in a) the cell $\omega^1$  and  in b) the cell $\omega^0$ .}%
\label{f1}%
\end{figure}

Composite materials are used extensively in the modern engineering practice.
They are mathematically interpreted, modelled and treated by applying the
Floquet-Bloch-Gelfand-theory (FBG-theory in the sequel)
for elliptic spectral boundary value problems
with periodic coefficients in periodic domains. This approach has led to many
important theoretical results and applications in topics like homogenization,
diffraction in waveguides, band-gap engineering etc. The usual setting for the
theory concerns purely periodic media (periodic coefficients, periodic
domains), though certain types of perturbations may be allowed as well. In
this paper we introduce and examine quite novel type of  perturbations of a
double-periodic medium with semi-infinite rows of foreign inclusions as
depicted in Fig.\,\ref{f1},\,a),\,b). The perturbation may influence
the essential spectrum of the problem, and the description of this spectrum
becomes the main goal of our paper. Such perturbations of periodic lattices
have not been the subject of  thorough mathematical investigation, but in  the physical
literature they   are related, for example, to defected photonic
crystals, cf. \cite[Ch.\,7]{book}.

In the case of homogeneous media, foreign inclusions
like infinite or semi-infinite strips are called
open waveguides; we accept the same terminology in the periodic case. The
physical meaning of the notion is clear, see \cite{Har1, Har2} for acoustics
and \cite{book, EKSS, HK1, HK2} for similar defects of periodic media in
solid state physics and optical systems. Our approach can readily be
generalized for different shapes of insertions, see Fig.\,\ref{f2}.
These types of perturbed periodic media may also appear in applications like
composite materials, because of improper manufacturing or also a specific
feature created on purpose. For technical simplicity we will deal with
particular open waveguides as depicted in Fig.\,\ref{f1}, but in
Section \ref{sect4} we describe minor modifications of our approach for shapes like in
Fig.\,\ref{f2}  as well as  other generalizations, including smoothness of coefficients and the boundary.

\begin{figure}
[ptb]
\begin{center}
\includegraphics[
height=1.4105in,
width=4.1148in
]%
{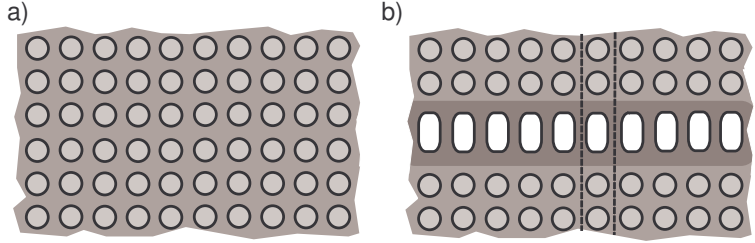}%
\caption{a) Purely periodic plane $\Omega^{0}$. b) Plane
$\Omega^{\sharp}$ with infinite perforated strip $\Omega^1$ (dark grey) and the strip
$\Pi$ (between dotted lines).}%
\label{f0}%
\end{center}
\end{figure}

We will consider a general elliptic spectral problem, see (\ref{15})--(\ref{16}), in
the perturbed domain $\Omega\subset{\mathbb{R}}^{2}$, Fig.\,\ref{f1}.
The main result of
the paper, Theorem \ref{thMAIN} in Section \ref{sect3.5}, contains
the following statement:
the essential spectrum
$\sigma_{es}\left(  \mathcal{A}\right)  $ of the main problem is the union of
two subsets, the first of which is the essential
spectrum $\sigma_{es}\left(  \mathcal{A}^{0}\right)  $ of the problem on the purely periodic domain $\Omega^0$, see Fig.\,\ref{f0},\,a), and  the second one, denoted by $\sigma^{\sharp} $, which is caused by  the perturbation of the domain, the open waveguide
$\Omega^{1+}$, see Fig.\,\ref{f1}. The subset  $\sigma^{\sharp} $ equals
the union of the discrete spectra of a family of model problems on a domain
$\Omega^\sharp$,  Fig.\,\ref{f0},\,b), which is related to $\Omega^{1+}$.
We will discuss the structure of $\sigma_{es}\left(  \mathcal{A}\right)  $
in more detail in Section \ref{sect1.5}, after presenting the notation and
definitions.

\noindent

\subsection{Purely periodic medium.\label{sect1.2}}

We consider the geometry of the unperturbed perforated plane and the
related spectral boundary value problem.
Let $\mathbf{Q}=\mathbf{Q}_{1}\subset{\mathbb{R}}^{2}$ be the unit open
square, $\mathbf{Q}_{a}=\left(  0,a\right)  \times\left(  0,a\right)  ,$ and
let $\vartheta^{0}$ be an open set in the plane $\mathbb{R}^{2}$ with closure
$\overline{\vartheta^{0}}=\vartheta^{0}\cup\partial\vartheta^{0}%
\subset\mathbf{Q}$ and a smooth boundary $\partial\vartheta$ of H\"{o}lder
class $C^{2,\delta}$, where $\delta>0$. Here $\vartheta^{0}$ is a not
necessarily connected set describing perforation; neither the case
$\vartheta^{0}=\emptyset$ is excluded. We define the periodicity cells
(outlined by dotted line in Fig.\,\ref{f1},\,b)
\begin{equation}
\omega^{0}\left(  \alpha\right)  =\left\{  x=\left(  x_{1},x_{2}\right)
:\left(  x_{1}-\alpha_{1},x_{2}-\alpha_{2}\right)  \in\omega^{0}%
:=\mathbf{Q}\setminus\overline{\vartheta^{0}}\right\}  \label{1}%
\end{equation}
where $\alpha=\left(  \alpha_{1},\alpha_{2}\right)  \in\mathbb{Z}^{2}$ is a
multi-index and $\mathbb{Z}=\left\{  0,\pm1,\pm2,...\right\}  .$ The
perforated plane $\Omega^{0}$ is covered by the periodicity cells (\ref{1}),
so it satisfies%
\begin{equation}
\overline{\Omega^{0}}=%
{\displaystyle\bigcup_{\alpha\in\mathbb{Z}^{2}}}
\overline{\omega^{0}\left(  \alpha\right)  },\ \ \ \ \Omega^{0}=\mathbb{R}%
^{2}\setminus%
{\displaystyle\bigcup_{\alpha\in\mathbb{Z}^{2}}}
\overline{\vartheta^{0}\left(  \alpha\right)  }, \label{2}%
\end{equation}
where $\vartheta^{0}\left(  \alpha\right)  $ denotes a translation of
$\vartheta^{0}$ similarly to (\ref{1}). Assuming $\omega^{0}$ to be a domain,
we observe that $\omega^{0}$ includes the boundary strip%
\begin{equation}
\left\{  x\in\mathbf{Q}\,:\,\mbox{dist}\left(  x,\partial\mathbf{Q}\right)
<d\right\}  \label{strip}%
\end{equation}
with some $d>0$ and, thus, $\Omega^{0}$ is a domain and in particular a
connected set as well.

We consider the spectral problem%
\begin{align}
L^{0}\left(  x,\nabla\right)  u^{0}\left(  x\right)   &  =\lambda u^{0}\left(
x\right)  ,\ \ \ x\in\Omega^{0},\label{3}\\
N^{0}\left(  x,\nabla\right)  u^{0}\left(  x\right)   &  =0,\ \ \ x\in
\partial\Omega^{0}, \label{4}%
\end{align}
and its variational form%
\begin{equation}
\left(  A^{0}D\left(  \nabla\right)  u^{0},D\left(  \nabla\right)
v^{0}\right)  _{\Omega^{0}}=\lambda\left(  u^{0},v^{0}\right)  _{\Omega^{0}%
},\ \ \ \forall v^{0}\in H^{1}\left(  \Omega^{0}\right)  ^{n}. \label{5}%
\end{equation}
Concerning the notation, by $\lambda$ we understand a spectral parameter
and $u=\left(  u_{1},...,u_{n}\right)  ^{\top}$ is a vector function, realized
as a column, so that $\top$ stands for transposition. The matrix differential
operators $L^{0}$ and $N^{0}$ in (\ref{3}) and (\ref{4}) take the form
\begin{align}
L^{0}\left(  x,\nabla\right)   &  =\overline{D\left(  -\nabla\right)  }^{\top
}A^{0}\left(  x\right)  D\left(  \nabla\right)  ,\label{6}\\
N^{0}\left(  x,\nabla\right)   &  =\overline{D\left(  \nu\left(  x\right)
\right)  }^{\top}A^{0}\left(  x\right)  D\left(  \nabla\right)  , \label{7}%
\end{align}
while $\nabla=\operatorname{grad},$ $\nu=\left(  \nu^{1},\nu^{2}\right)
^{\top}$ is the unit outward normal vector on $\partial\Omega^{0}$ and $A^{0}$
is a matrix function of size $m\times m$ with $C^{1,\delta}$-smooth entries;
$A^{0}\left(  x\right)  $ is Hermitian, positive definite and 1-periodic in
$x\in\mathbb{R}^{2}$ (it is convenient to define the coefficient matrix in the
intact plane). Furthermore, $D\left(x \right)  $ is a  matrix function of size $m\times n$ and
linear in $x=(x_{1}, x_{2})$, and $D(0,0)
=\mathbb{O}_{n\times m}$ is the null matrix. The substitution $x_{j}%
\mapsto \partial/\partial x_{j}$ gives an $m\times
n$-matrix $D\left(  \nabla\right)  $ of first order differential operators with constant complex
coefficients, while $L^{0}(x,\nabla)$ (respectively,  $N^{0}(x,\nabla)$)  is
a  formally self-adjoint differential matrix operator in divergence form (resp. the
corresponding Neumann condition operator); the bar in (\ref{6})
and (\ref{7}) indicates complex conjugation. Finally, $\left(  \ ,\ \right)
_{\Omega^{0}}$ is the natural scalar product in the Lebesgue space
$L^{2}\left(  \Omega^{0}\right)  $, and $H^{1}\left(  \Omega^{0}\right)  $ is
the Sobolev space with standard norm%
\[
\left\Vert w;H^{1}\left(  \Omega^{0}\right)  \right\Vert =\left(  \left\Vert
\nabla w;L^{2}\left(  \Omega^{0}\right)  \right\Vert ^{2}+\left\Vert
w;L^{2}\left(  \Omega^{0}\right)  \right\Vert ^{2}\right)  ^{1/2}.
\]
The last superscript $n$ in the integral identity (\ref{5}) shows the number
of components in the test function $v=\left(  v_{1},...,v_{n}\right)  ^{\top
},$ however, this index is omitted in the notation of norms and scalar
products. In this way, the left- and right-hand sides of (\ref{5}) involve
scalar products in $L^{2}\left(  \Omega^{0}\right)  ^{m}$ and $L^{2}\left(
\Omega^{0}\right)  ^{n}$, respectively.

We assume that  $D\left(x\right)  $ 
is \textit{algebraically complete}
\cite{Nec}: there exists a $\varrho_{D}\in\mathbb{N}=\left\{
1,2,3,...\right\}  $ such that for any row $p=\left(  p_{1},...,p_{n}\right)
$ of homogeneous polynomials in $\xi=\left(  \xi_{1},\xi_{2}\right)
\in\mathbb{R}^{2}$ of common degree $\varrho\geq\varrho_{D},$ one can find a
row of polynomials $q=\left(  q_{1},...,q_{n}\right)  $ satisfying
\begin{equation}
p\left(  \xi\right)  =q\left(  \xi\right)  D\left(  \xi\right)
,\ \ \ \ \forall\xi\in\mathbb{R}^{2}. \label{8}%
\end{equation}
According to \cite[\S \,3.7.4]{Nec}, this assumption
yields the Korn inequality%
\begin{equation}
\left\Vert u;H^{1}\left(  \Omega^{0}\right)  \right\Vert ^{2}\leq
c_{D,\omega^{0}}\left(  \left\Vert D\left(  \nabla\right)  u;L^{2}\left(
\omega^{0}\right)  \right\Vert ^{2}+\left\Vert u;L^{2}\left(  \omega
^{0}\right)  \right\Vert ^{2}\right)  \label{9}%
\end{equation}
and, hence, the sesquilinear Hermitian positive form $a^{0}\left(  u^{0}%
,v^{0}\right)  $ on the left in (\ref{5}) is closed in $H^{1}\left(
\Omega^{0}\right)  ^{n}.$ Moreover, the operator $L^{0}$ is elliptic and the
Neumann boundary operator $N^{0}$ covers it in the Shapiro-Lopatinskii sense
everywhere on $\partial\Omega^{0}$ (see, e.g., \cite[Thm.\,1.9]{na262}).

Owing to the above-mentioned properties of $a^{0}$, the problem (\ref{3}),
(\ref{4}) is associated with a positive self-adjoint operator $\mathcal{A}%
^{0}$ in $L^{2}\left(  \Omega^{0}\right)  ^{n}$ with the differential
expression $L^{0}\left(  x,\nabla\right)  $ and the domain%
\begin{equation}
\mathcal{D}\left(  \mathcal{A}^{0}\right)  =\left\{  u^{0}\in H^{2}\left(
\Omega^{0}\right)  :N^{0}\left(  x,\nabla\right)  u^{0}\left(  x\right)
=0,\ x\in\partial\Omega^{0}\right\}  ; \label{10}%
\end{equation}
see \cite[Ch.\,10]{BiSo} and \cite[Ch.\,13]{Ru}. The description of the spectrum
$\sigma\left(  \mathcal{A}^{0}\right)  $ is well known and will be presented
in Section \ref{sect2}.

We emphasize that the results of the paper remain valid for other types
of boundary conditions, in particular, for the Dirichlet conditions, cf. Section
\ref{sect4.1}.$(i)$,$(iv)$. A description of all admissible boundary conditions
can be found in \cite{LiMa} and \cite[\S\,1]{na262}. Moreover, the
$C^{2,\delta}$-smoothness of the boundary was assumed in order to have the
elementary formula (\ref{10}) for the operator domain $\mathcal{D}\left(
\mathcal{A}^{0}\right)  $ and to simplify the technical computations.
The results of our paper  hold true in the case of uniformly Lipschitz boundaries,
but the proofs would  require small modifications, cf. Section \ref{sect4.2}.

\subsection{Periodic medium with open semi-infinite periodic
waveguide.\label{sect1.3}}

In this section we describe the geometry of the open waveguide and the full
spectral problem.
Let $\varpi^{+}$ be the semi-strip  $\left\{  x:x_{1}>0,\ \left\vert x_{2}\right\vert <h\right\}  $  with some
$h\in\mathbb{N}$ (overshaded in Fig.\,\ref{f1},\,b). In the rectangle $\mathbf{Q}_{1}^{h}=\left(  0,1\right)
\times\left(  -h,h\right)  ,$ we introduce an open set $\vartheta^{1}$ with a
smooth boundary $\partial\vartheta^{1}$ and closure $\overline{\vartheta^{1}%
}\subset\mathbf{Q}_{1}^{h}.$ We define a semi-infinite row of holes or
inclusions as depicted in Fig.\,\ref{f1},\,a) or b), respectively:%
\begin{align}
\omega^{1}\left(  \alpha_{1}\right)   &  =\left\{  x:\left(  x_{1}-\alpha
_{1},x_{2}\right)  \in\omega^{1}\right\}  ,\ \ \alpha_{1}\in\mathbb{N}%
_{0}=\mathbb{N\cup}\left\{  0\right\}  ,\label{11}\\
\omega^{1}  &  =\mathbf{Q}_{1}^{h}\setminus\overline{\vartheta^{1}%
},\ \ \vartheta^{1}\left(  \alpha_{1}\right)  =\left\{  x:\left(  x_{1}%
-\alpha_{1},x_{2}\right)  \in\vartheta^{1}\right\}  .\nonumber
\end{align}
The cell $\omega^{1}$ is outlined in Fig.\,\ref{f1},\,a), by dotted line. We also
introduce a smooth Hermitian matrix function $A^{1+}$ in $\mathbb{R}^{2},$
entries of which are supported in $\varpi^{+}$ and become $1$-periodic in
$x_{1}$ inside the semi-strip $\left\{  x\in\varpi^{+}:x_{1}>R\right\}  ,$
$R\in\mathbb{N}$,
\begin{equation}
A^{1+}\left(  x\right)  =A^{1}\left(  x\right)  \text{ \ for \ }%
x_{1}>R,\ \ \ A^{1}\left(  x_{1}\pm1,x_{2}\right)  =A^{1}\left(  x_{1}%
,x_{2}\right)  . \label{A}%
\end{equation}
Furthermore,
\begin{equation}
\Omega^{1+}=\varpi^{+}\setminus%
{\displaystyle\bigcup_{\alpha_{1}\in\mathbb{N}_{0}}}
\overline{\vartheta^{1}\left(  \alpha_{1}\right)  },\ \ \ \overline
{\Omega^{1+}}=%
{\displaystyle\bigcup_{\alpha_{1}\in\mathbb{N}_{0}}}
\overline{\omega^{1}\left(  \alpha_{1}\right)  }, \label{12}%
\end{equation}
and the sum%
\begin{equation}
A\left(  x\right)  =A^{0}\left(  x\right)  +A^{1+}\left(  x\right)  \label{13}%
\end{equation}
is assumed to be positive definite in $\overline{\Omega}$, where%
\begin{equation}
\Omega=\left(  \Omega^{0}\setminus\varpi^{+}\right)  \cup\Omega^{1+}.
\label{14}%
\end{equation}
In other words, we make a perturbation of coefficients and boundary inside the
semi-strip $\varpi^{+}$, see Fig.\,\ref{f1}. For instance, one may suppose that
$A^{1+}$ is the null matrix $\mathbb{O}_{m\times m}$ and $\vartheta^{0}%
\neq\varnothing,$ $\vartheta^{1}=\varnothing,$ $h=1,$ which means filling in
all holes inside $\varpi^{+},$ cf. Fig.\,\ref{f1},\,b). Vice versa, in the case
$\vartheta^{0}=\varnothing,$ $\vartheta^{1}\neq\varnothing,$ one perforates
the plane $\mathbb{R}^{2}$ with a semi-infinite row of holes, see
Fig.\,\ref{f1},\,a). Even in the case of absence of holes we still call (\ref{12}) the
perforated strip, and the perforated plane $\Omega^{0}$ in (\ref{2}) can also
contain no holes.

Notice that the Korn inequality%
\begin{equation}
\left\Vert u;H^{1}\left(  \Omega\right)  \right\Vert ^{2}\leq c_{D}\left(
\omega^{0},\omega^{1}\right)  \left(  \left\Vert D\left(  \nabla\right)
u;L^{2}\left(  \Omega\right)  \right\Vert ^{2}+\left\Vert u;L^{2}\left(
\Omega\right)  \right\Vert ^{2}\right)  \label{AK}%
\end{equation}
is still valid and can be derived by summing up inequalities of type (\ref{9})
in the cells $\omega^{0}\left(  \alpha\right)  $ with $\alpha\in\mathbb{Z}%
^{2}$ ($\alpha_{1}>0,$ $-h\leq\alpha_{2}<h$ excluded) and $\omega^{1}\left(
\alpha_{1}\right)  $ with $\alpha_{1}\in\mathbb{N}_{0}.$

Replacing $A^{0}\left(  x\right)  $ with (\ref{13}) in (\ref{6}) still gives
an elliptic and formally self-adjoint matrix operator $L\left(  x,\nabla
_{x}\right)  .$ The same change in (\ref{7}) yields the Neumann boundary
condition operator $N\left(  x,\nabla_{x}\right)  ,$ where $\nu$ is regarded
as the outward unit normal on $\partial\Omega$. In the domain (\ref{14}) we
consider the spectral problem
\begin{align}
L\left(  x,\nabla\right)  u\left(  x\right)   &  =\lambda u\left(  x\right)
,\ \ \ x\in\Omega,\label{15}\\
N\left(  x,\nabla\right)  u\left(  x\right)   &  =0,\ \ \ x\in\partial\Omega,
\label{16}%
\end{align}
and the corresponding integral identity
\begin{equation}
a\left(  u,v\right)  :=\left(  AD\left(  \nabla\right)  u,D\left(
\nabla\right)  v\right)  _{\Omega}=\lambda\left(  u,v\right)  _{\Omega
},\ \ \ \forall v\in H^{1}\left(  \Omega\right)  ^{n}. \label{17}%
\end{equation}
Since $a$ remains as a closed positive Hermitian form in $H^{1}\left(
\Omega\right)  ^{n},$ the variational formulation (\ref{17}) of the problem
(\ref{15}), (\ref{16}) supplies it with a positive self-adjoint operator
$\mathcal{A}$ in $L^{2}\left(  \Omega\right)  ^{n}$ with the differential
expression $L\left(  x,\nabla\right)  $ and the domain%
\begin{equation}
\mathcal{D}\left(  \mathcal{A}\right)  =\left\{  u\in H^{2}\left(  \Omega
^{0}\right)  ^{n}:N\left(  x,\nabla\right)  u\left(  x\right)  =0,\ x\in
\partial\Omega\right\}  ; \label{18}%
\end{equation}
see again \cite[Ch.\,10]{BiSo} and \cite[Ch.\,13]{Ru}. This notation is quite
similar to the one used in Section \ref{sect1.2}.

The main goal of the paper is to describe the essential component $\sigma
_{es}\left(  \mathcal{A}\right)  $ in the spectrum of $\mathcal{A}$,
\begin{equation}
\sigma\left(  \mathcal{A}\right)  =\sigma_{di}\left(  \mathcal{A}\right)
\cup\sigma_{es}\left(  \mathcal{A}\right)  . \label{19}%
\end{equation}
We emphasize that in general
\begin{equation}
\sigma_{es}\left(  \mathcal{A}\right)  \neq\sigma_{es}\left(  \mathcal{A}%
^{0}\right)  , \label{spnot}%
\end{equation}
and moreover, the discrete spectrum of $\mathcal{A}^{0}$ is empty, thus,
\begin{equation}
\sigma_{es}\left(  \mathcal{A}^{0}\right)  =\sigma\left(  \mathcal{A}%
^{0}\right)  . \label{spA}%
\end{equation}
We will identify the difference%
\begin{equation}
\sigma_{ad}\left(  \mathcal{A}\right)  =\sigma_{es}\left(  \mathcal{A}\right)
\setminus\sigma\left(  \mathcal{A}^{0}\right)  , \label{spdif}%
\end{equation}
but leave aside two interesting and important questions. First,
we are not able to describe  completely the component $\sigma_{di}\left(
\mathcal{A}\right)  $ in a general perturbed problem (\ref{15}), (\ref{16}),
although, of course concrete examples of isolated and embedded eigenvalues in
$\sigma_{po}\left(\mathcal{A}\right)$ can be constructed in scalar problems. Second, the
existence or absence of the point spectrum $\sigma_{po}\left(  \mathcal{A}%
^{0}\right)  $ (eigenvalues of infinite multiplicity) in the purely periodic
problem (\ref{3}), (\ref{4}) remains unknown; note that this question is
answered in the literature only for particular scalar problems (see, e.g.,
papers \cite{Fil1, Fil2, Fil3} and books \cite{ReSi, Kuchbook}). Notice that
$\sigma_{po}\left(  \mathcal{A}^{0}\right)  $ can be included in $\sigma
_{po}\left(  \mathcal{A}\right)  $, but the latter stays unknown, too.

\subsection{Discussion on the main result\label{sect1.5}}

In the case of the purely periodic plane $\Omega^0$, Fig.\,\ref{f0},\,a),
the spectrum $\sigma_{es}\left(  \mathcal{A}^{0}\right)  $ of the problem (\ref{3})-(\ref{4})
has  representation as a union of spectral bands (see (\ref{33}), (\ref{34}), below),
which    is a well-known consequence of the FBG-theory; we refer here to \cite{KuchUMN,
Kuchbook, Skri}.  Consider for a moment the domain $\Omega^\sharp$ of
Fig.\,\ref{f0},\,b) with foreign inclusions or holes, which form a periodic row, infinite in both directions. Then, the spectrum may be different from the purely periodic case,
and we denote by $\sigma^\sharp$ the essential spectrum of the problem
on $\Omega^\sharp$ (cf. (\ref{54}), below).
Analogously to \cite{EKSS, HK1, HK2}, the
increment $\sigma^\sharp \setminus  \sigma_{es}\left(  \mathcal{A}^{0}\right) $  
can be detected by performing the partial FBG-transform in $x_{1}$-direction
(Section \ref{sect3.1}) and investigating the kernel of the model problem  in the
perforated strip $\Pi$ (separated by dashed lines in Fig.\,\ref{f0},\,b); cf. (\ref{59}), (\ref{58})\,).
This problem depends on the Floquet parameter $\zeta$, and, for certain
values of the spectral parameter $\lambda\in\left(  0,+\infty\right)  $ and $\zeta\in\left[  0,2\pi\right),$
it can have a solution in the Sobolev space $H^{2}(\Pi)^{n}.$
Such values  $\lambda$ form the point spectrum $\sigma_{po}\left(
\mathcal{A}^{\sharp}\left(  \zeta\right)  \right)  $ of the model problem \eqref{88}
for the operator $\mathcal{A}^{\sharp}\left(  \zeta\right)  $.

Our main result in Theorem \ref{thMAIN} says that the following formula
holds true for the problem (\ref{15}), (\ref{16}) in the
periodic plane with the immersed $\mathsf{I}$-shaped\ open waveguide, see
(\ref{14}) and Fig.\,\ref{f1},\,b):
\begin{equation}
\sigma_{ad}\left(  \mathcal{A}\right)  =
\sigma^\sharp \setminus  \sigma_{es}\left(  \mathcal{A}^{0}\right)
= {\displaystyle\bigcup_{\zeta\in\left[  0,2\pi\right)  }}
\sigma_{po}\left(  \mathcal{A}^{\sharp}\left(  \zeta\right)  \right)
\setminus \sigma_{es}\left(  \mathcal{A}^{0}\right)  .\label{sigma}%
\end{equation}
The last set in this formula requires some comments.
First, embedded eigenvalues in $\sigma_{em}\left(  \mathcal{A}^{\sharp}\left(
\zeta\right)  \right)  =\sigma_{po}\left(  \mathcal{A}^{\sharp}\left(
\zeta\right)  \right)  \setminus\sigma_{di}\left(  \mathcal{A}^{\sharp}\left(
\zeta\right)  \right)  $ live inside the essential spectrum $\sigma
_{es}\left(  \mathcal{A}^{\sharp}\left(  \zeta\right)  \right)  $, which
in turn is contained in $\sigma_{es}\left(  \mathcal{A}^{0}\right)  $ (compare (\ref{90})
with (\ref{33})). Second, $\sigma_{es}\left(  \mathcal{A}^{\sharp}\left(
\zeta\right)  \right)  $ depends on $\zeta$ and therefore some points of the
discrete spectrum $\sigma_{di}\left(  \mathcal{A}^{\sharp}\left(
\zeta\right)  \right)  $ may fall into $\sigma_{es}\left(  \mathcal{A}%
^{\sharp}\left(  \zeta^{\prime}\right)  \right)  $ with $\zeta^{\prime}%
\neq\zeta.$ In any case none of the indicated points in $\sigma_{po}\left(
\mathcal{A}^{\sharp}\left(  \zeta\right)  \right)  $ stays in the
increment component (\ref{spdif}) of $\sigma_{es}\left(  \mathcal{A}\right)
.$

We emphasize that, contrary to the case of $\Omega^\sharp$,
the lacking periodicity of the domain $\Omega$ prevents a direct
use of the partial FBG-transform, hence, the proof of
\eqref{sigma} requires improved  mathematical tools.
The new procedure of our  paper involves the construction
of a parametrix and singular Weyl sequences in order to describe, respectively, the regularity
field and the essential spectrum of the operator $\mathcal{A}$.

\subsection{Cranked and branching open waveguides\label{sect1.4}}

\begin{figure}[ptb]
\begin{center}
\includegraphics[
height=3in,
width=3.6in
]{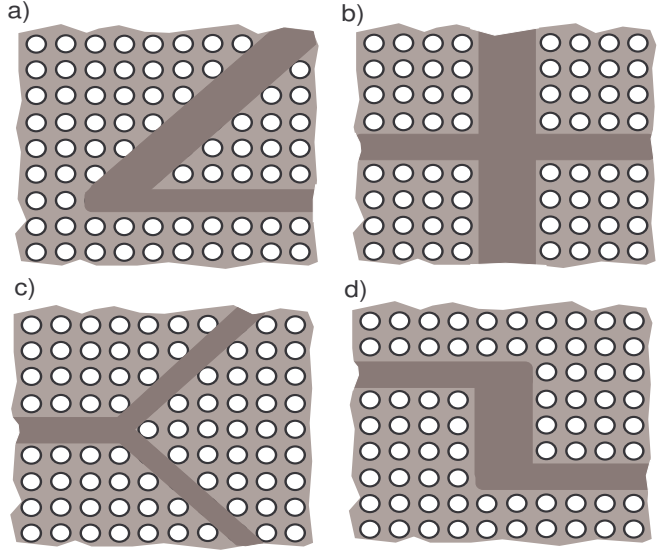}
\end{center}
\caption{ $\mathsf{V, X}$ and $\mathsf{Y, Z}$-shaped joints of open
waveguides (dark grey).}%
\label{f2}%
\end{figure}

Fig.\,\ref{f2} shows open waveguides of the shape of the letters $\mathsf{V},$
$\mathsf{X}$ and $\mathsf{Y},$ $\mathsf{Z}.$ They appear due to perforation
and perturbation of coefficients in overshadowed joints of semi-strips. In
Section \ref{sect4} we explain how the results of Section \ref{sect3} for the
$\mathsf{I}$-shaped case, Fig.\,\ref{f1},\,a),\,b), can be readily adapted to
these cranked and branching open waveguides. Here, skewed branches of the
$\mathsf{V}$- and $\mathsf{Y}$-shaped waveguides must maintain the periodicity and
thus the tangent of tilt angles has to be rational number. We do not know a
formula for the essential spectrum in the irrational case.

We also emphasize that for a clear reason, no relevant perturbation in a disk
$\mathbb{B}_{R}=\left\{  x:r=\left\vert x\right\vert <R\right\}  $ with radius
$R>0$ can affect the essential spectrum of the boundary value problem
(\ref{15}), (\ref{16}). Moreover, assume that $\widetilde{A}\left(  x\right)
$ is a matrix decaying together with its derivatives at infinity as $\left(
1+\left\vert x\right\vert \right)  ^{-\delta},$ $\delta>0$, and that the
coefficient matrix (\ref{13})%
\begin{equation}
A\left(  x\right)  +\widetilde{A}\left(  x\right)  \label{20}%
\end{equation}
still keeps the above mentioned basic properties of $A$. This replacement of
the coefficient matrix does not change $\sigma_{es}$.

All these generalizations and some others will be commented in Section
\ref{sect4}. We have chosen the very particular open waveguide in Section
\ref{sect1.3} in order to simplify the presentation, to illuminate the main
points of our approach and to avoid unimportant but cumbersome technical details.

\subsection{Structure of the paper\label{sect1.6}}

In Section \ref{sect2} we recall generally known information on the purely
periodic case which will be used later on. The main interest is focused
on the model problem (\ref{29}) in the periodicity cell $\varpi$, which is
obtained using the FBG-transform \cite{Gel}. The open periodic semi-infinite
waveguide will be  considered in Section \ref{sect3}, where we apply the partial
$(FGB)$-transform to formulate another model problem (\ref{88}) in the perforated
infinite strip $\Pi$ with periodicity conditions on its lateral sides.

The spectra of those two model problems form the essential spectrum
$\sigma_{es}\left(  \mathcal{A}\right)  $ of the original problem (\ref{17}).
To verify the corresponding formulas (\ref{A1}) and (\ref{34}), (\ref{A2}) we
first present two types of singular Weyl sequences for
$\mathcal{A}$ and on the other hand construct a right parametrix for the
formally self-adjoint problem (\ref{15}), (\ref{16}). This is the most
involved part of our paper. To that end, we follow \cite{na17} and also
\cite[\S 3.4]{NaPl}, and study an operator family for a second model problem
in the weighted Sobolev spaces (Kondratiev spaces) $W_{\beta}^{l}\left(
\Pi\right)  ,$ leading to important conclusions on
exponential decay properties of the solutions in the strip $\Pi.$ Finally, we
glue the parametrix (\ref{A9}) from solutions of the model problems with the
help of appropriate cut-off functions. The parametrix enables to prove that,
for any $\lambda\in\overline{\mathbb{R}_{+}}$ outside the union of sets
(\ref{34}) and (\ref{A2}), the operator of the inhomogeneous problem
(\ref{15}), (\ref{16}), cf. (\ref{A3}), is Fredholm in the
Sobolev-Slobodetskii spaces; therefore such points $\lambda$ form the
intersection of the regularity field of ${\mathcal{A}}$ with the semi-axis
$\overline{\mathbb{R}_{+}}.$ This completes the proof of Theorem \ref{thMAIN},
the central assertion in the paper.

We start the last section of the paper by describing several concrete problems
in acoustics, elasticity and piezoelectricity, to which our theory may apply.
However, as mentioned above,
the original exact formulation of the spectral problem was simplified in
several aspects, so we comment in the next subsections on certain supplementary
issues in order to obtain more generality for further interesting physical
applications.
We finish the paper with Section \ref{sect4.5}, where we present small
modifications of the parametrix to be applied to semi-bounded open periodic
waveguides in the shape of the letters $\mathsf{V,X}$ and $\mathsf{Y,Z}$ as in
Fig.\,\ref{f2}.

\section{Spectrum of the purely periodic problem\label{sect2}}

\subsection{Floquet-Bloch-Gelfand-transform\label{sect2.1}}

The spectrum of the operator $\mathcal{A}^0$ in the purely periodic
domain $\Omega^0$ can be studied with the help of the FBG-transform. This
will lead to the formula \eqref{34}, the main object of Section \ref{sect2}.
The FBG-transform is defined by %
\begin{equation}
u\mapsto\widehat{u}\left(  x;\eta\right)  =\frac{1}{2\pi}\sum_{\alpha
\in\mathbb{Z}^{2}}e^{-i\eta\cdot\left(  x+\alpha\right)  }u\left(
x+\alpha\right)  ,\ \ \ x\in\omega^{0}, \label{21}%
\end{equation}
and it establishes the isometric isomorphism%
\begin{equation}
L^{2}\left(  \Omega^{0}\right)  \cong L^{2}\left(  \mathbf{Y};L^{2}\left(
\omega^{0}\right)  \right)  \label{22}%
\end{equation}
(cf. \cite{Gel} and, e.g., \cite{Skri, Kuchbook}), where $\eta\cdot x=\eta
_{1}x_{1}+\eta_{2}x_{2},$%
\begin{equation}
\mathbf{Y}=\left\{  \eta=\left(  \eta_{1},\eta_{2}\right)  :\eta_{1},\eta
_{2}\in\left[  0,2\pi\right)  \right\}  ,\ \ \ \overline{\mathbf{Y}}%
=\overline{\mathbf{Q}_{2\pi}} , \label{YY}%
\end{equation}
and $L^{2}\left(  \mathbf{Y};\mathfrak{B}\right)  $ is the Lebesgue space of
abstract functions in $\eta\in\mathbf{Y}$ with values in Banach space
$\mathfrak{B}$\ and the norm%
\begin{equation}
\left\Vert U;L^{2}\left(  \mathbf{Y};\mathfrak{B}\right)  \right\Vert =\left(
\int_{0}^{2\pi}\int_{0}^{2\pi}\left\Vert U\left(  \eta\right)  ;\mathfrak{B}%
\right\Vert ^{2}d\eta_{1}d\eta_{2}\right)  ^{1/2}. \label{23}%
\end{equation}
Moreover, the mapping
\begin{equation}
H^{2}\left(  \Omega^{0}\right)  ^{n}\ni u\mapsto\widehat{u}\in L^{2}\left(
\mathbf{Y};H_{\mathrm{per}}^{2}\left(  \omega^{0}\right)  ^{n}\right)
\label{24}%
\end{equation}
is an isomorphism, too. Here, $H_{\mathrm{per}}^{2}\left(  \omega^{0}\right)
^{n}$ is the subspace of vector functions $U\in H^{2}\left(  \omega
^{0}\right)  ^{n}$ satisfying the periodicity conditions
\begin{align}
U\left(  0,x_{2}\right)   &  =U\left(  1,x_{2}\right)  ,\ \ \ \frac{\partial
U}{\partial x_{1}}\left(  0,x_{2}\right)  =\frac{\partial U}{\partial x_{1}%
}\left(  1,x_{2}\right)  ,\ \ \ x_{2}\in\left(  0,1\right)  ,\label{25}\\
U\left(  x_{1},0\right)   &  =U\left(  x_{1},1\right)  ,\ \ \ \frac{\partial
U}{\partial x_{2}}\left(  x_{1},0\right)  =\frac{\partial U}{\partial x_{2}%
}\left(  x_{1},1\right)  ,\ \ \ x_{1}\in\left(  0,1\right)  .\nonumber
\end{align}
In what follows we shorten the notation $\partial U/\partial x_{j}$ to
$\partial_{j}U.$

The inverse FBG-transform is given by%
\begin{align}
U  &  =\widehat{u}\mapsto u\left(  x\right)  =\frac{1}{2\pi}\int_{0}^{2\pi
}\int_{0}^{2\pi}e^{i\eta\cdot x}U\left(  x-\left[  x\right]  ;\eta\right)
d\eta_{1}d\eta_{2},\label{26}\\
\left[  x\right]   &  =\left(  \left[  x_{1}\right]  ,\left[  x_{2}\right]
\right)  ,\ \ \ \left[  t\right]  =\max\left\{  \tau\in\mathbb{Z}:\tau\leq
t\right\}  .\nonumber
\end{align}

\subsection{The model problem on the periodicity cell\label{sect2.2}}

Owing to (\ref{21}), we have%
\begin{equation}
\widehat{\,\,\nabla^{\alpha}v\,\,}\left(  x;\eta\right)  =\left(  \partial
_{1}+i\eta\right)  ^{\alpha_{1}}\left(  \partial_{2}+i\eta\right)
^{\alpha_{2}}\widehat{v}\left(  x;\eta\right)  ,\ \ \ \nabla^{\alpha}%
=\partial_{1}^{\alpha_{1}}\partial_{2}^{\alpha_{2}} \label{der}%
\end{equation}
and thus the FBG-transform (\ref{21}) converts the problem (\ref{3}),
(\ref{4}) into the following problem, depending on the parameter $\eta
\in\mathbf{Y}$, in the periodicity cell $\omega^{0}$ ,%
\begin{align}
L^{0}\left(  x,\nabla+i\eta\right)  U\left(  x;\eta\right)   &  =\Lambda
\left(  \eta\right)  U\left(  x;\eta\right)  ,\ \ \ x\in\omega^{0}%
,\label{27}\\
N^{0}\left(  x,\nabla+i\eta\right)  U\left(  x;\eta\right)   &  =0,\ \ \ x\in
\gamma^{0}=\partial\omega^{0}\cap\mathbf{Q}, \label{28}%
\end{align}
together with the periodicity conditions (\ref{25}) on the exterior part
$\partial\mathbf{Q}$ of the boundary of the cell $\omega^{0}.$ The variational
formulation of problem (\ref{27}), (\ref{28}), (\ref{25}) amounts to finding a
number $\Lambda\left(  \eta\right)  \in\mathbb{C}$ and a non-trivial vector
function $U\left(  \cdot;\eta\right)  \in H_{\mathrm{per}}^{1}\left(
\omega^{0}\right)  ^{n}$ such that
\begin{equation}
\left(  A^{0}D\left(  \nabla+i\eta\right)  U,D\left(  \nabla+i\eta\right)
V\right)  _{\omega^{0}}=\Lambda\left(  \eta\right)  \left(  U,V\right)
_{\omega^{0}},\ \ \ \forall V\in H_{\mathrm{per}}^{1}\left(  \omega
^{0}\right)  ^{n}. \label{29}%
\end{equation}
In view of the compact embedding $H^{1}\left(  \omega^{0}\right)  \subset
L^{2}\left(  \omega^{0}\right)  $ in the bounded domain $\omega^{0},$ the
spectrum $\sigma_{\eta}^{0}$ of the variational problem (\ref{29}) and
boundary value problem (\ref{27}), (\ref{28}), (\ref{25}) is discrete and
forms the unbounded monotone sequence%
\begin{equation}
0\leq\Lambda_{1}\left(  \eta\right)  \leq\Lambda_{2}\left(  \eta\right)
\leq...\leq\Lambda_{k}\left(  \eta\right)  \leq...\rightarrow+\infty,
\label{30}%
\end{equation}
where multiplicities are counted. According to the general results of the
perturbation theory for linear operators\footnote{A quadratic pencil easily
reduces to a linear non-self-adjoint spectral family.}, the functions%
\begin{equation}
\mathbf{Y}\ni\eta\mapsto\Lambda_{k}\left(  \eta\right)  \label{31}%
\end{equation}
are continuous, see for example \cite{HiPhi, Kato}. Moreover, they are $2\pi
$-periodic in $\eta_{1}$ and $\eta_{2}$, because for any eigenpair $\left\{
\Lambda\left(  \eta\right)  ,U\left(  x;\eta\right)  \right\}  $ of the
problem (\ref{27}), (\ref{28}), (\ref{25}) at some $\eta\in\mathbf{Y},$
\begin{equation}
\left\{  \Lambda\left(  \eta\right)  ,e^{\pm2\pi ix_{p}}U\left(
x;\eta\right)  \right\}  ,\ p=1,2, \label{32}%
\end{equation}
remains an eigenpair of the same problem but at $\eta\pm2\pi e_{\left(
p\right)  }$ where $e_{\left(  p\right)  }=\left(  \delta_{1,p},\delta
_{2,p}\right)  ^{\top}\in\mathbb{R}^{2}$ is the unit vector of the $\eta_{p}$-axis.

The above mentioned properties of functions (\ref{31}) ensure that the
\textit{spectral bands}
\begin{equation}
B_{k}^{0}=\left\{  \Lambda_{k}\left(  \eta\right)  :\eta\in\mathbf{Y}\right\}
,\ k\in\mathbb{N}, \label{33}%
\end{equation}
are bounded connected closed segments in $\overline{\mathbb{R}_{+}}.$ The
formula%
\begin{equation}
\sigma\left(  \mathcal{A}^{0}\right)  =
\sigma_{es}\left(  \mathcal{A}^{0}\right)  =
{\displaystyle\bigcup_{k\in\mathbb{N}}}
B_{k}^{0} \label{34}%
\end{equation}
for the spectrum of the problem (\ref{3}), (\ref{4}) is well-known, see, e.g.,
\cite{KuchUMN, Skri, Kuchbook}, but we briefly comment on its proof in
Sections \ref{sect2.3}--\ref{sect2.4}, since we will need some of these arguments later.

\subsection{Unique solution of the inhomogeneous problem\label{sect2.3}}
As regards to \eqref{34}, we now prove the inclusion
$\sigma\left(  \mathcal{A}^{0}\right)  \subset
\cup_{k\in\mathbb{N}} B_{k}^{0}$.
Let us consider the boundary value problem%
\begin{align}
L^{0}\left(  x,\nabla\right)  u^{0}\left(  x\right)  -\lambda u^{0}\left(
x\right)   &  =f^{0}\left(  x\right)  ,\ \ \ x\in\Omega^{0},\label{35}\\
N^{0}\left(  x,\nabla\right)  u^{0}\left(  x\right)   &  =g^{0}\left(
x\right)  ,\ \ \ x\in\partial\Omega^{0},\nonumber
\end{align}
with the data%
\begin{equation}
f^{0}\in L^{2}\left(  \Omega^{0}\right)  ^{n},\ g^{0}\in H^{1/2}\left(
\partial\Omega^{0}\right)  ^{n} \label{36}%
\end{equation}
and a fixed parameter $\lambda=\lambda^{0}$ such that%
\begin{equation}
\lambda^{0}\notin B_{k}^{0},\ \ \forall k\in\mathbb{N}. \label{37}%
\end{equation}
In (\ref{36}), $H^{1/2}\left(  \partial\Omega^{0}\right)  $ stands for the
Sobolev-Slobodetskii space of traces with the intrinsic norm%
\begin{equation}
\left\Vert g^{0};H^{1/2}\left(  \partial\Omega^{0}\right)  \right\Vert
=\inf\left\{  \left\Vert G^{0};H^{1}\left(  \Omega^{0}\right)  \right\Vert
:G^{0}=g^{0}\text{ on }\partial\Omega^{0}\right\}  . \label{38}%
\end{equation}
This norm is equivalent to the following one:%
\[
\left(  \left\Vert g^{0};L^{2}\left(  \partial\Omega^{0}\right)  \right\Vert
^{2}+\int_{\partial\Omega^{0}}\int_{\partial\Omega_{l}^{0}\left(  x\right)
}\left\vert g^{0}\left(  x\right)  -g^{0}\left(  y\right)  \right\vert
^{2}\frac{ds_{x}ds_{y}}{\left\vert x-y\right\vert ^{2}}\right)  ^{1/2}.
\]
Here, $ds_{x}$ is the arc length element on $\partial\Omega^{0}$ and
$\partial\Omega_{l}^{0}\left(  x\right)  =\left\{  y\in\partial\Omega
^{0}:\left\vert y-x\right\vert <l\right\}  $ while $l\in\left(  0,+\infty
\right)  $ can be fixed arbitrarily.

Clearly, the mapping%
\begin{equation}
H^{2}\left(  \Omega^{0}\right)  ^{n}\ni u\mapsto\left\{  L^{0}u^{0}%
-\lambda^{0}u^{0},N^{0}u^{0}\right\}  \in L^{2}\left(  \Omega^{0}\right)
^{n}\times H^{1/2}\left(  \partial\Omega^{0}\right)  ^{n} \label{39}%
\end{equation}
is continuous for any $\lambda^{0}$, but in the case (\ref{37}) it becomes an
isomorphism. Indeed, the FBG-transform (\ref{21}) turns (\ref{35}) into the
parameter-dependent problem%
\begin{align}
L^{0}\left(  x,\nabla+i\eta\right)  {\widehat{u}\,}^{0}\left(  x;\eta\right)
-\lambda^{0}{\widehat{u}\,}^{0}\left(  x;\eta\right)   &  ={\widehat{f}\ }%
^{0}\left(  x;\eta\right)  ,\ \ \ x\in\omega^{0},\label{40}\\
N^{0}\left(  x,\nabla+i\eta\right)  {\widehat{u}\,}^{0}\left(  x;\eta\right)
&  ={\widehat{g}\;}^{0}\left(  x;\eta\right)  ,\ \ \ x\in\gamma^{0},\nonumber
\end{align}
with the periodicity conditions (\ref{25}). The right-hand sides meet the
estimate%
\begin{align}
&  \left\Vert {\widehat{f}\ }^{0};L^{2}\left(  0,2\pi;L^{2}\left(  \omega
^{0}\right)  \right)  \right\Vert ^{2}+\left\Vert {\widehat{g}\;}^{0}%
;L^{2}\left(  0,2\pi;H^{1/2}\left(  \gamma^{0}\right)  \right)  \right\Vert
^{2}\label{41}\\
&  \leq c\left(  \left\Vert f^{0};L^{2}\left(  \Omega^{0}\right)  \right\Vert
^{2}+\left\Vert g^{0};H^{1/2}\left(  \partial\Omega^{0}\right)  \right\Vert
^{2}\right)  ,\nonumber
\end{align}
while the necessary information about the Sobolev-Slobodetskii space is
provided by the isomorphisms (\ref{22}), (\ref{23}) and the formula
(\ref{der}) for derivatives. By the assumption (\ref{37}), the problem
(\ref{40}), (\ref{25}) has for any $\eta\in\mathbf{Y}$ a unique solution
denoted by
\begin{equation}
{\widehat{u}\,}^{0}\left(  x;\eta\right)  =R^{0}\left(  \lambda^{0}%
;\eta\right)  \left\{  {\widehat{f}\ }^{0}\left(  \cdot;\eta\right)
,{\widehat{g}\;}^{0}\left(  \cdot;\eta\right)  \right\}  \label{42}%
\end{equation}
and estimated as follows:%
\begin{equation}
\left\Vert {\widehat{u}\,}^{0}\left(  \cdot;\eta\right)  ;H^{2}\left(
\omega^{0}\right)  \right\Vert ^{2}\leq C\left(  \left\Vert {\widehat{f}%
\ }^{0}\left(  \cdot;\eta\right)  ;L^{2}\left(  \omega^{0}\right)  \right\Vert
^{2}+\left\Vert {\widehat{g}\;}^{0}\left(  \cdot;\eta\right)  ;H^{1/2}\left(
\gamma^{0}\right)  \right\Vert ^{2}\right)  . \label{43}%
\end{equation}
Since the constant $C$ does not depend on $\eta\in\mathbf{Y}$, it suffices to
apply the inverse FBG-transform (\ref{26}) and to derive from (\ref{43}) and
(\ref{41}) the inequality%
\begin{equation}
\left\Vert u^{0};H^{2}\left(  \Omega^{0}\right)  \right\Vert ^{2}\leq c\left(
\left\Vert f^{0};L^{2}\left(  \Omega^{0}\right)  \right\Vert ^{2}+\left\Vert
g^{0};H^{1/2}\left(  \partial\Omega^{0}\right)  \right\Vert ^{2}\right)
\label{44}%
\end{equation}
for the unique solution $u^{0}\in H^{2}\left(  \Omega^{0}\right)  ^{n}$ of the
problem (\ref{35}) with fixed parameter (\ref{37}). Thus, mapping (\ref{39})
with this $\lambda^{0}$ is indeed an isomorphism, which in particular means
that $\lambda^{0}$ belongs to the regularity field of the operator
$\mathcal{A}^{0}$ in (\ref{10}). This coincides with the resolvent set of
$\mathcal{A}^{0}$ because the discrete spectrum $\sigma_{di}\left(
\mathcal{A}^{0}\right)  $ is evidently empty.

\begin{remark}
\textrm{\label{remSoSl}Dealing with the Sobolev-Slobodetskii norms in
(\ref{41}), (\ref{43}) and (\ref{44}), it is much more convenient to use the
definition (\ref{38}) and some extensions $G$ of $g$ and $\widehat{G}\left(
\cdot;\eta\right)  $ of $\widehat{g}\left(  \cdot;\eta\right)  $ such that%
\begin{align*}
\left\Vert g;H^{1/2}\left(  \partial\Omega^{0}\right)  \right\Vert  &
\leq2\left\Vert G;H^{1}\left(  \Omega^{0}\right)  \right\Vert \text{,}\\
\left\Vert \widehat{g};L^{2}\left(  0,2\pi;H^{1/2}\left(  \gamma^{0}\right)
\right)  \right\Vert  &  \leq2\left\Vert \widehat{G};L^{2}\left(
0,2\pi;H^{1/2}\left(  \gamma^{0}\right)  \right)  \right\Vert
.\ \ \blacksquare
\end{align*}
}
\end{remark}

\subsection{The singular Weyl sequence\label{sect2.4}}

We next show that
$\sigma\left(  \mathcal{A}^{0}\right)  \supset
\cup_{k\in\mathbb{N}} B_{k}^{0}$.
Let us assume%
\begin{equation}
\lambda^{0}\in B_{k}^{0}\ \text{for some }k\in\mathbb{N} \label{45}%
\end{equation}
so that there exist $\eta^{0}\in\mathbf{Y}$ and $U^{0}\left(  \cdot;\eta
^{0}\right)  \in H^{2}\left(  \omega^{0}\right)  $ such that $\left\{
\eta^{0},U^{0}\left(  \cdot;\eta^{0}\right)  \right\}  $ is an eigenpair of
the problem (\ref{27}), (\ref{28}), (\ref{25}) with $\Lambda\left(  \eta
^{0}\right)  =\lambda^{0}$. By a direct calculation, one easily deduces that
the Bloch wave%
\begin{equation}
u^{0}\left(  x\right)  =e^{i\eta^{0}\cdot x}U^{0}\left(  x;\eta^{0}\right)
\label{46}%
\end{equation}
satisfies the differential equations (\ref{3}) and the boundary conditions
(\ref{4}), although it does of course not fall into the Sobolev space
$H^{2}\left(  \Omega^{0}\right)  ^{n}.$ However, (\ref{46}) is useful for
constructing a singular sequence $\left\{  u^{0j}\right\}  _{j\in\mathbb{N}}$
in $\mathcal{D}\left(  \mathcal{A}^{0}\right)  \subset H^{2}\left(  \Omega
^{0}\right)  ^{n}$ for the operator $\mathcal{A}^{0}$ at the point
$\lambda^{0},$ namely a sequence with the following properties:

1$^{\text{o}}$ \ $\left\Vert u^{0j};L^{2}\left(  \Omega^{0}\right)
\right\Vert =1;$

2$^{\text{o}}$ \ $u^{0j}\rightharpoondown0$ weakly in $L^{2}\left(  \Omega
^{0}\right)  ^{n}$ as $j\rightarrow+\infty;$

3$^{\text{o}}$ \ $\left\Vert \mathcal{A}^{0}u^{0j}-\lambda^{0}u^{0j}%
;L^{2}\left(  \Omega^{0}\right)  \right\Vert \rightarrow0$ as $j\rightarrow
+\infty.$

\bigskip

To define the entries of this sequence, we introduce the plateau function%
\begin{equation}
\chi_{j}\left(  t\right)  =\chi\left(  t-2^{j}\right)  \chi\left(
2^{j+1}-t\right)  , \label{47}%
\end{equation}
where $\chi\in C^{\infty}\left(  \mathbb{R}\right)  $ is a cut-off function
such that%
\begin{equation}
\chi\left(  t\right)  =1\text{ for\ }t\geq d,\ \ \ \ \ \ \chi\left(  t\right)
=0\text{ \ for \ }t\leq0,\text{ \ \ \ \ \ }0\leq\chi\leq1 \label{48}%
\end{equation}
and $d>0$ is taken from (\ref{strip}); therefore, in the vicinity of each
component $\gamma^{0}\left(  \alpha\right)  =\partial\vartheta^{0}\left(
\alpha\right)  $ of the boundary $\partial\Omega^{0}$ the two-dimensional
plateau function%
\begin{equation}
X_{j}\left(  x\right)  =\chi_{j}\left(  x_{1}\right)  \chi_{j}\left(
x_{2}\right)  \label{49}%
\end{equation}
becomes a constant, either $1$ or $0.$

We set%
\begin{equation}
v^{0j}\left(  x\right)  =X_{j}\left(  x\right)  u^{0}\left(  x\right)
,\ \ \ u^{0j}\left(  x\right)  =\left\Vert v^{0j};L^{2}\left(  \Omega
^{0}\right)  \right\Vert ^{-1}v^{0j}\left(  x\right)  . \label{50}%
\end{equation}
The above specification of $X_{j}$ shows that%
\[
N^{0}\left(  x,\nabla\right)  v^{0j}\left(  x\right)  =X_{j}\left(  x\right)
N^{0}\left(  x,\nabla\right)  u^{0}\left(  x\right)  =0,\ \ \ x\in
\partial\Omega^{0},
\]
and, hence, $u^{0j}\in\mathcal{D}\left(  \mathcal{A}^{0}\right)  .$

\begin{figure}[ptb]
\begin{center}
\includegraphics[
height=2in,
width=5in
]{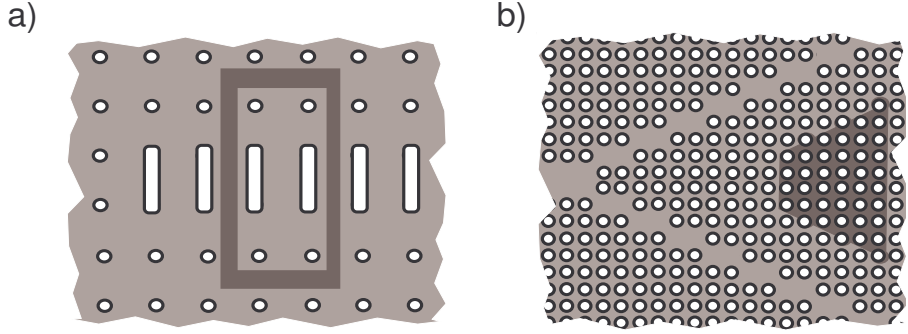}
\end{center}
\caption{Supports  of the cut-off functions.}%
\label{f4}%
\end{figure}

The property 1$^{\text{o}}$ clearly holds true. Furthermore, the weak
convergence to $0$ in 2$^{\text{o}}$ occurs at least along a subsequence of
indices $j\in\mathbb{N}$, because, by (\ref{47}) and (\ref{49}),
supp$\,\left(  v^{0j}\right)  \cap\mbox{supp}\left(  v^{0l}\right)
=\varnothing$ as $j\neq l.$ It remains to verify the property 3$^{\text{o}}.$
Recalling (\ref{50}) and (\ref{46}), we have%
\begin{align}
\left\Vert v^{0j};L^{2}\left(  \Omega^{0}\right)  \right\Vert ^{2}  &
\geq\sum_{\alpha\in\mathbb{Z}^{2}\cap\left[  2^{j}+1,2^{j+1}-1\right]  }%
\int_{\omega^{0}\left(  \alpha\right)  }\left\vert U^{0}\left(  x;\eta
^{0}\right)  \right\vert ^{2}dx=\label{51}\\
&  =\left(  2^{j+1}-2^{j}-2\right)  ^{2}\left\Vert U^{0}\left(  \cdot;\eta
^{0}\right)  ;L^{2}\left(  \omega^{0}\right)  \right\Vert ^{2}\geq c_{0}%
2^{2j}\nonumber
\end{align}
with $c_{0}>0$ and $j\geq2.$ At the same time, we obtain%
\[
f^{0j}=L^{0}v^{0j}-\lambda^{0}v^{0j}=\left[  L^{0},X_{j}\right]  u^{0}%
+X_{j}\left(  L^{0}u^{0}-\lambda^{0}u^{0}\right)  =\left[  L^{0},X_{j}\right]
u^{0}%
\]
where $\left[  L^{0},X_{j}\right]  $ stands for the commutator of the
differential operator (\ref{6}) and the multiplication operator with $X_{j}.$
Owing to definition (\ref{47})-(\ref{49}), the plateau function (\ref{49})
varies only inside the union of four rectangles of size $d\times2^{j}$ and the
common area $O\left(  2^{j}\right)  $ (see the overshaded frame in
Fig.\,\ref{f4},\,a)\,). Hence, due to periodicity in (\ref{46}), we arrive at the
inequality%
\[
\left\Vert f^{0j};L^{2}\left(  \Omega^{0}\right)  \right\Vert ^{2}\leq
c_{1}2^{j},
\]
which together with (\ref{51}) and (\ref{50}) prove the relation%
\[
\left\Vert \mathcal{A}^{0}u^{0j}-\lambda^{0}u^{0j};L^{2}\left(  \Omega
^{0}\right)  \right\Vert \leq C2^{-j/2},\ \ C=c_{1}/c_{0}%
\]
as well as the property 3$^{\text{o}}$.

By the Weyl criterion (cf. \cite[Th.\,9.12]{BiSo} or \cite[Th.\,VII.12]%
{ReSi1}), the point (\ref{45}) lives in the essential spectrum $\sigma
_{es} \left( \mathcal{A}^{0} \right).$ This and the material of Section \ref{sect2.3} confirm the formula
(\ref{34}).

\section{Spectrum of the open waveguide in periodic medium\label{sect3}}

\subsection{Partial Floquet-Bloch-Gelfand-transform\label{sect3.1}}
Our aim is to apply the partial FBG-transform to detect the effect of the
open waveguide to the essential spectrum of the problem \eqref{15}--\eqref{16}.
Due to the lack of periodicity, this cannot be done directly in $\Omega$,
hence, we introduce and study in Sections \ref{sect3.1}--\ref{sect3.3}
the problem in the domain $\Omega^\sharp$, Fig.\,\ref{f0},\,b).
The results of Section \ref{sect3.3} will be applied in Sections
\ref{sect3.4}--\ref{sect3.5} to the original problem, which leads to the
proof of the main result, Theorem \ref{thMAIN}.

Similarly to (\ref{12})-(\ref{14}) we introduce the infinite periodically
perforated strip $\Omega^{1}$ (overshaded in Fig.\,\ref{f0},\,b),%
\begin{equation}
\Omega^{1}=\varpi\setminus%
{\displaystyle\bigcup_{\alpha_{1}\in\mathbb{Z}}}
\overline{\omega^{1}\left(  \alpha_{1}\right)  },\ \ \ \ \varpi=\mathbb{R}%
\times\left(  -h,h\right) ,  \label{opi}%
\end{equation}
and the positive definite Hermitian matrix%
\begin{equation}
A^{\sharp}\left(  x\right)  =A^{0}\left(  x\right)  +A^{1}\left(  x\right)
\label{52}%
\end{equation}
which happens to be 1-periodic in the variable $x_{1}.$ In the domain%
\begin{equation}
\Omega^{\sharp}=\left(  \Omega^{0}\setminus\varpi\right)  \cup\Omega^{1}
\label{53}%
\end{equation}
we consider the auxiliary boundary-value problem%
\begin{align}
L^{\sharp}\left(  x,\nabla\right)  u^{\sharp}\left(  x\right)  -\lambda
^{\sharp}u^{\sharp}\left(  x\right)   &  =f^{\sharp}\left(  x\right)
,\ \ \ x\in\Omega^{\sharp},\label{54}\\
N^{\sharp}\left(  x,\nabla\right)  u^{\sharp}\left(  x\right)   &  =g^{\sharp
}\left(  x\right)  ,\ \ \ x\in\partial\Omega^{\sharp},\nonumber
\end{align}
where $\lambda^{\sharp}$ is a fixed parameter and the operators $L^{\sharp
},N^{\sharp}$ are given by (\ref{6}), (\ref{7}) with the change $A^{0}\mapsto
A^{\sharp}.$

The domain (\ref{53}) is also 1-periodic along the $x_{1}$-axis, that is%
\begin{equation}
\Omega^{\sharp}=\left\{  x:\left(  x_{1}\pm1,x_{2}\right)  \in\Omega^{\sharp
}\right\}  , \label{sharp}%
\end{equation}
and we define the perforated strip bounded by dashed lines in
Fig.\,\ref{f0},\,b),
\begin{equation}
\Pi=\left\{  x\in\Omega^{\sharp}:x_{1}\in\left(  0,1\right)  \right\}  .
\label{55}%
\end{equation}

The periodicity observed in (\ref{52}) and (\ref{53}) allows us to apply the
partial FBG-transform, see e.g. \cite{na17},%
\begin{equation}
u^{\sharp}\mapsto U^{\sharp}\left(  x;\zeta\right)  =\frac{1}{\sqrt{2\pi}}%
\sum_{\alpha_{1}\in\mathbb{Z}}e^{-i\zeta \left(  x_{1}+\alpha_{1}\right)
}u^{\sharp}\left(  x_{1}+\alpha_{1};x_{2}\right)  \label{56}%
\end{equation}
which establishes the isomorphisms%
\begin{equation}
L^{2}\left(  \Omega^{\sharp}\right)  \simeq L^{2}\left(  0,2\pi;L^{2}\left(
\Pi\right)  \right)  ,\ \ \ H^{2}\left(  \Omega^{\sharp}\right)  \approx
L^{2}\left(  0,2\pi;H_{\mathrm{per}\sharp}^{2}\left(  \Pi\right)  \right)
\label{57}%
\end{equation}
where the first one is isometric (cf. (\ref{22})), while $H_{\mathrm{per}%
\sharp}^{2}\left(  \Pi\right)  ^{n}$ is the subspace of functions $U^{\sharp
}\in H^{2}\left(  \Pi\right)  ^{n}$ satisfying the periodicity conditions on
the lateral sides of the perforated strip \eqref{55}
\begin{equation}
U^{\sharp}\left(  0;x_{2}\right)  =U^{\sharp}\left(  1;x_{2}\right)
,\ \ \ \partial_{1}U^{\sharp}\left(  0;x_{2}\right)  =\partial_{1}U^{\sharp
}\left(  1;x_{2}\right)  ,\ \ \ x_{2}\in\mathbb{R}. \label{58}%
\end{equation}
The inverse partial FBG-transform is given by%
\begin{equation}
u^{\sharp}\left(  x_{1},x_{2}\right)  =\frac{1}{\sqrt{2\pi}}\int_{0}^{2\pi
}e^{i\zeta x_{1}}U^{\sharp}\left(  x_{1}-\left[  x_{1}\right]  ,x_{2}%
;\zeta\right)  d\zeta, \label{Gin}%
\end{equation}
see \cite{na17} and cf. (\ref{21}), (\ref{26}).

The FBG-transform (\ref{56}) applies to the problem (\ref{54}) and
converts it into the parameter-dependent problem
\begin{align}
L^{\sharp}\left(  x,\partial_{1}+i\zeta,\partial_{2}\right)  U^{\sharp}\left(
x;\eta\right)  -\lambda^{\sharp}U^{\sharp}\left(  x;\eta\right)   &
=F^{\sharp}\left(  x;\eta\right)  ,\ \ \ x\in\Pi,\label{59}\\
N^{\sharp}\left(  x,\partial_{1}+i\zeta,\partial_{2}\right)  U^{\sharp}\left(
x;\eta\right)   &  =G^{\sharp}\left(  x;\eta\right)  ,\ \ \ x\in
\Gamma,\nonumber
\end{align}
with the periodicity conditions (\ref{58}). Here, $\lambda^{\sharp}%
\in\mathbb{R}_{+}$ is fixed, $F^{\sharp},G^{\sharp}$ are the FBG-images of
$f^{\sharp},g^{\sharp}$ and%
\begin{equation}
\Gamma=\left\{  x\in\partial\Pi:0<x_{1}<1\right\}  \label{60}%
\end{equation}
is the interior boundary of the perforated strip $\Pi$.

\subsection{Second model problem in weighted spaces\label{sect3.2}}

To study the problem (\ref{59}), (\ref{58}), we introduce the weighted Sobolev
space $W_{\beta}^{2}\left(  \Pi\right)  $ (the exponential Kondratiev space
\cite{Ko}) as a completion of the linear set $C_{c}^{\infty}\left(
\overline{\Pi}\right)  $ (infinitely differentiable functions with compact
supports) with respect to the norm%
\begin{equation}
\left\Vert v;W_{\beta}^{2}\left(  \Pi\right)  \right\Vert =\left(  \left\Vert
\nabla^{2}v;L_{\beta}^{2}\left(  \Pi\right)  \right\Vert ^{2}+\left\Vert
\nabla v;L_{\beta}^{2}\left(  \Pi\right)  \right\Vert ^{2}+\left\Vert
v;L_{\beta}^{2}\left(  \Pi\right)  \right\Vert ^{2}\right)  ^{1/2}, \label{71}%
\end{equation}
where $\nabla^{2}v$ is the family of all second-order derivatives of $v,$
$\beta\in\mathbb{R}$ is a weight index and $L_{\beta}^{2}\left(  \Pi\right)  $
stands for the weighted Lebesgue space,%
\begin{equation}
\left\Vert v;L_{\beta}^{2}\left(  \Pi\right)  \right\Vert =\left\Vert
e^{\beta\left\vert x_{2}\right\vert }v;L^{2}\left(  \Pi\right)  \right\Vert .
\label{72}%
\end{equation}
Notice that $W_{\beta}^{2}\left(  \Pi\right)  $ consists of all functions
$v\in H_{loc}^{2}\left(  \overline{\Pi}\right)  $ with finite norm (\ref{71}).
The norm modified by omitting $\left\Vert \nabla v;L_{\beta}^{2}\left(
\Pi\right)  \right\Vert $ on the right in (\ref{71}) remains equivalent to the
original one. For $\beta=0,$ we have $W_{0}^{2}\left(  \Pi\right)
=H^{2}\left(  \Pi\right)  $, but in the case $\beta>0$ ($\beta<0)$ the
Kondratiev space includes functions with an exponential decay (growth) at
infinity with decay (growth) rate controlled by the weight index. By
$W_{\beta,per\sharp}^{2}\left(  \Pi\right)  ,$ we understand the subspace of
functions subject to the periodicity conditions (\ref{58}), and $W_{\beta
}^{1/2}\left(  \Gamma\right)  $ is the weighted Sobolev-Slobodetskii space
with the intrinsic norm%
\begin{align}
\Big\Vert v;W_{\beta}^{1/2}(\Gamma)\Big\Vert& =\inf  \Big\{\left\Vert
V;W_{\beta}^{1}\left(  \Pi\right)  \right\Vert
\nonumber\\
& =\Big(\left\Vert \nabla
V;L_{\beta}^{2}\left(  \Pi\right)  \right\Vert ^{2}
+\left\Vert V;L_{\beta}^{2}\left(  \Pi\right)  \right\Vert
^{2}\Big)^{1/2}\,:\,V=v\ \mathrm{on}\ \Gamma\Big\}. \label{73}%
\end{align}
We emphasize that $W_{\beta}^{1/2}\left(  \Gamma\right)  $ does not require
periodicity conditions because (\ref{60}) includes the interior part of the boundary.

Owing to definitions (\ref{71})-(\ref{73}) and formula (\ref{57}), the partial
FBG-transform establishes the isomorphisms%
\begin{align}
L_{\beta}^{2}\left(  \Omega^{\sharp}\right)   &  \simeq L^{2}\left(
0,2\pi;L_{\beta}^{2}\left(  \Pi\right)  \right)  ,\ \ \ W_{\beta}^{2}\left(
\Omega^{\sharp}\right)  \approx L^{2}\left(  0,2\pi;W_{\beta,per\sharp}%
^{2}\left(  \Pi\right)  \right)  ,\label{74}\\
W_{\beta}^{3/2}\left(  \partial\Omega^{\sharp}\right)   &  \approx
L^{2}\left(  0,2\pi;W_{\beta}^{3/2}\left(  \Gamma^{\sharp}\right)  \right)
.\nonumber
\end{align}
The problem operator for (\ref{59}), (\ref{58}) ,
\begin{align}
W_{\beta,per\sharp}^{2}\left(  \Pi\right)  ^{n}  &  \ni U^{\sharp}\mapsto
T_{\beta}^{\sharp}\left(  \lambda^{\sharp},\zeta\right)  U^{\sharp}\nonumber\\
&  =\left\{  L^{\sharp}\left(  x,\partial_{1}+i\zeta,\partial_{2}\right)
-\lambda^{\sharp},N^{\sharp}\left(  x,\partial_{1}+i\zeta,\partial_{2}\right)
\right\}  \in L_{\beta}^{2}\left(  \Pi\right)  ^{n}\times W_{\beta}%
^{1/2}\left(  \Gamma\right)  ^{n} \label{75}%
\end{align}
is continuous for any $\beta\in\mathbb{R}$ and $\zeta\in\mathbb{C}$. However,
it has better properties under additional assumptions described in terms of
the operator pencil%
\begin{gather}
\mathbb{C}\ni\eta_{2}\mapsto\mathfrak{A}_{\Lambda}\left(  \zeta,\eta
_{2}\right)  =\left\{  L^{0}\left(  x,\partial_{1}+i\zeta,\partial_{2}%
+i\eta_{2}\right)  -\Lambda,N^{0}\left(  x,\partial_{1}+i\zeta,\partial
_{2}+i\eta_{2}\right)  \right\}  :\label{76}\\
H_{\mathrm{per}}^{2}\left(  \omega^{0}\right)  ^{n}\rightarrow L^{2}\left(
\omega^{0}\right)  ^{n}\times H^{1/2}\left(  \gamma^{0}\right)  ^{n}\nonumber
\end{gather}
which corresponds to the problem (\ref{27}), (\ref{28}), (\ref{25}) with the
fixed parameter $\Lambda=\lambda^{\sharp}\in\mathbb{R}$ and the dual
FBG-variable%
\begin{equation}
\eta=\left(  \zeta,\eta_{2}\right)  \in\mathbf{Y}, \label{77}%
\end{equation}
where $\zeta\in\left[  0,2\pi\right)  $ is taken from (\ref{59}). Concerning
$\eta_{2}$ as a spectral parameter, we regard (\ref{76}) as a quadratic pencil
in $\eta_{2},$ which is a particular case of a holomorphic spectral family,
see \cite[Ch.1]{GoKr}.

The following assertion is proved in \cite{na17}, see also \cite[Thms 3.4.7
and 5.1.4]{NaPl}.

\begin{theorem}
\label{thna17}The operator (\ref{75}) is Fredholm if and only if the segment%
\begin{equation}
\Upsilon_{\beta}=\left\{  \eta_{2}\in\mathbb{C}:\operatorname{Re}\eta_{2}%
\in\left[  0,2\pi\right)  ,\ \operatorname{Im}\eta_{2}=\beta\right\}
\label{78}%
\end{equation}
in the complex plane $\mathbb{C}$ is free of the spectrum of the pencil
(\ref{76}). If a point of the spectrum belongs to $\Upsilon_{\beta},$ then the
range of the operator $T_{\beta}^{\sharp}\left(  \lambda^{\sharp}%
,\zeta\right)  $ is not closed.
\end{theorem}

Let us assume that $\lambda^{\sharp}=\lambda^{0}$ satisfies (\ref{37}) and, in
particular, that the segment $\Upsilon_{0}$ is free of the spectrum of
$\mathfrak{A}_{\lambda^{\sharp}}\left(  \zeta,\cdot\right)  $ for all
$\zeta\in\left[  0,2\pi\right)  .$ Indeed, if $\eta_{2}\in\Upsilon_{0}$
belongs to the spectrum, then $\lambda^{\sharp}=\Lambda\left(  \zeta,\eta
_{2}\right)  $ becomes an eigenvalue of the problem (\ref{27}), (\ref{28}),
(\ref{25}) and therefore falls into some spectral band $B_{k}.$ By the
analytic Fredholm alternative, see, e.g., \cite[Thm.\,1.5.1]{GoKr} or
\cite[Thm.\,VI.14]{ReSi1}, we conclude that the spectrum of the pencil
(\ref{76}) consists of a countable set of (normal) eigenvalues without finite
accumulation points. For $\eta_{2}\in\mathbb{C}$, the spectrum is invariant
with respect to the shifts $\eta_{2}\mapsto\eta_{2}\pm2\pi$ along the real
axis, by the same argument as in (\ref{31}) and (\ref{32}). Moreover, it is
mirror symmetric with respect to the real axis because, for fixed real
$\lambda^{\sharp}$ and $\zeta,$ the problems (\ref{27}), (\ref{28}),
(\ref{25}) with $\eta_{2}$ and $\overline{\eta_{2}}$ are formally adjoint.
Finally, under the assumption (\ref{37}), there exists a positive continuous
$2\pi$-periodic function%
\begin{equation}
\left[  0,2\pi\right)  \ni\zeta\mapsto\beta_{0}\left(  \lambda^{\sharp}%
;\zeta\right)  \in\mathbb{R}_{+} \label{80}%
\end{equation}
such that the rectangle%
\begin{equation}
\left\{  \eta_{2}\in\mathbb{C}:\operatorname{Re}\eta_{2}\in\left[
0,2\pi\right)  ,\ \left\vert \operatorname{Im}\eta_{2}\right\vert <\beta
_{0}\left(  \lambda^{\sharp};\zeta\right)  \right\}  \label{81}%
\end{equation}
does not contain any eigenvalue of the pencil $\mathfrak{A}_{\lambda^{\sharp}%
}\left(  \zeta,\cdot\right)  $, but the segments $\Upsilon_{\pm\beta
_{0}\left(  \lambda^{\sharp};\zeta\right)  }$ surely do. We further put%
\begin{equation}
\beta_{0}\left(  \lambda^{\sharp}\right)  =\min\left\{  \beta_{0}\left(
\lambda^{\sharp};\zeta\right)  :\zeta\in\left[  0,2\pi\right)  \right\}  >0.
\label{82}%
\end{equation}

When $\lambda^{\sharp}$ and $\zeta$ are real, the problem (\ref{59}),
(\ref{58}) in the infinite strip $\Pi$ is formally self-adjoint, which just
means the validity of the Green formula%
\begin{align}
&  \left(  L^{\sharp}\left(  x,\partial_{1}+i\zeta,\partial_{2}\right)
U-\lambda^{\sharp}U,V\right)  _{\Pi}+\left(  N^{\sharp}\left(  x,\partial
_{1}+i\zeta,\partial_{2}\right)  U,V\right)  _{\Gamma}\label{79}\\
&  =\left(  U,L^{\sharp}\left(  x,\partial_{1}+i\zeta,\partial_{2}\right)
V-\lambda^{\sharp}V\right)  _{\Pi}+\left(  U,N^{\sharp}\left(  x,\partial
_{1}+i\zeta,\partial_{2}\right)  V\right)  _{\Gamma}\nonumber
\end{align}
for all $U,V\in H_{\mathrm{per}\sharp}^{1}\left(  \Pi\right)  ^{n}$. Hence,
our assumption (\ref{37}) and Theorem \ref{thna17} assure that the operator
$T_{0}\left(  \lambda^{\sharp};\zeta\right)  $ in the Sobolev space
$H_{\mathrm{per}\sharp}^{1}\left(  \Pi\right)  ^{n}$ is Fredholm of index
zero. The next theorem follows from a general result in \cite{na17} (see also
\cite[\S 3.4 and \S 5.1]{NaPl}) and it concerns the finite-dimensional
subspace%
\begin{align}
\ker T_{\beta}\left(  \lambda^{\sharp};\zeta\right)  =  &  \big\{U^{\sharp}\in
W_{\beta,per\sharp}^{1}\left(  \Pi\right)  ^{n}:U^{\sharp}\text{ satisfies the
homogeneous}\label{83}\\
&  \ \text{problem}\ (\ref{59}),(\ref{58})\ \text{with}\ F^{\sharp
}=0,G^{\sharp}=0\big\}.\nonumber
\end{align}

\begin{theorem}
\label{thna17.1} Under the condition (\ref{37}) with $\lambda^{\sharp}%
=\lambda^{0}$, the subspace $\ker T_{\beta}\left(  \lambda^{\sharp}%
;\zeta\right)  $ is independent of the weight index $\beta\in\left(
-\beta_{0}\left(  \lambda^{\sharp};\zeta\right)  ,\beta_{0}\left(
\lambda^{\sharp};\zeta\right)  \right)  $, where $\beta_{0}\left(
\lambda^{\sharp};\zeta\right)  $ is determined in (\ref{80})-(\ref{81}).
\end{theorem}

\begin{theorem}
\label{thna17.2} Let $\lambda^{\sharp}$ satisfy (\ref{37}) and let $\beta
\in\left(  -\beta_{0}\left(  \lambda^{\sharp};\zeta\right)  ,\beta_{0}\left(
\lambda^{\sharp};\zeta\right)  \right)  $ where $\beta_{0}\left(
\lambda^{\sharp};\zeta\right)  $ is taken from (\ref{80}). For any fixed
$\zeta\in\left[  0,2\pi\right)  $, the problem (\ref{59}), (\ref{58}) with the
right-hand side%
\begin{equation}
\left\{  F^{\sharp},G^{\sharp}\right\}  \in L_{\beta}^{2}\left(  \Pi\right)
^{n}\times W_{\beta}^{1/2}\left(  \Gamma\right)  ^{n} \label{84}%
\end{equation}
has a solution $U^{\sharp}\in W_{\beta,per\sharp}^{2}\left(  \Pi\right)  ^{n}%
$, if and only if the compatibility conditions%
\begin{equation}
\left(  F^{\sharp},V\right)  _{\Pi}+\left(  G^{\sharp},V\right)  _{\Gamma
}=0,\ \ \ \forall V\in\ker T_{-\beta}\left(  \lambda^{\sharp};\zeta\right)
\label{85}%
\end{equation}
is met. This solution is defined up to an addendum in $\ker T_{0}\left(
\lambda^{\sharp};\zeta\right)  =\ker T_{\pm\beta}\left(  \lambda^{\sharp
};\zeta\right)  $ and, if it is subject to the orthogonality condition%
\begin{equation}
\left(  U^{\sharp},V\right)  _{\Pi}=0,\ \ \ \forall V\in\ker T_{0}\left(
\lambda^{\sharp};\zeta\right)  , \label{86}%
\end{equation}
then it becomes unique and meets the estimate%
\begin{equation}
\left\Vert U^{\sharp};W_{\beta}^{2}\left(  \Pi\right)  \right\Vert \leq
c\left(  \beta\right)  \left(  \left\Vert F^{\sharp};L_{\beta}^{2}\left(
\Pi\right)  \right\Vert +\left\Vert G^{\sharp};W_{\beta}^{1/2}\left(
\Gamma\right)  \right\Vert \right)  . \label{87}%
\end{equation}
In the case $|\beta|<\beta_{0}(\lambda^{\natural})$ the constant $c(\beta)$
can be chosen independently of $\zeta\in\lbrack0,2\pi)$.
\end{theorem}

These theorems mean that any solution $U^{\sharp}\in H_{\mathrm{per}\sharp
}^{2}\left(  \Pi\right)  ^{n}$ of the homogeneous problem (\ref{59}),
(\ref{58}) has exponential decay at infinity. Moreover, such solutions
form the set of all defect functionals (\ref{85}) of the operator $T_{\beta
}\left(  \lambda^{\sharp};\zeta\right)  $ with $\left\vert \beta\right\vert
<\beta_{0}\left(  \lambda^{\sharp};\zeta\right)  .$ Finally, assume that
$U_{\left(  0\right)  }^{\sharp}\in H_{\mathrm{per}\sharp}^{2}\left(
\Pi\right)  ^{n}$ and $U_{\left(  \beta\right)  }^{\sharp}\in W_{\beta
^{\sharp},per\sharp}^{2}\left(  \Pi\right)  ^{n},$ $\beta^{\sharp}\in\left(
0,\beta_{0}\left(  \lambda^{\sharp};\zeta\right)  \right)  ,$ are solutions of
the problem (\ref{59}),(\ref{58}), and the right-hand side $\left\{
F^{\sharp},G^{\sharp}\right\}  $ satisfies (\ref{84}) for both $\beta=0$ and
$\beta=\beta^{\sharp}.$ Then, clearly, $U_{\left(  0\right)  }^{\sharp}$ and
$U_{\left(  \beta\right)  }^{\sharp}$ may differ by an element of subspace
(\ref{83}) only, and if the orthogonality condition (\ref{86}) holds true for
both, then $U_{\left(  0\right)  }^{\sharp}=U_{\left(  \beta\right)  }%
^{\sharp}$.

\subsection{The spectrum of the model problem in the strip\label{sect3.3}}

Let us consider the spectral problem in $\Pi,$ which is the homogeneous
($F^{\sharp}=0,$ $G^{\sharp}=0$) problem (\ref{59}), (\ref{58}) for the
spectral parameter $\lambda^{\sharp}.$ Its variational formulation is: find
a number $\lambda^{\sharp}\in\mathbb{C}$ and a non-trivial vector function
$U^{\sharp}\in H_{\mathrm{per}\sharp}^{1}\left(  \Pi\right)  ^{n}$ such that%
\begin{equation}
\left(  A^{\sharp}D\left(  \partial_{1}+i\zeta,\partial_{2}\right)  U^{\sharp
},D\left(  \partial_{1}+i\zeta,\partial_{2}\right)  V^{\sharp}\right)  _{\Pi
}=\lambda^{\sharp}\left(  U^{\sharp},V^{\sharp}\right)  _{\Pi},\ \ \ \forall
V^{\sharp}\in H_{\mathrm{per}\sharp}^{1}\left(  \Pi\right)  ^{n}. \label{88}%
\end{equation}

\begin{figure}[ptb]
\begin{center}
\includegraphics[
height=1.9718in,
width=1.3785in
]{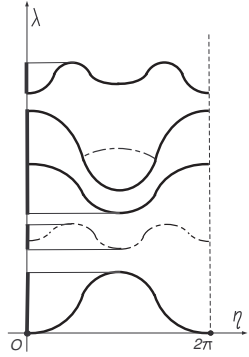}
\end{center}
\caption{Graphs of eigenvalues of the problems in $\varpi^{0}$ and $\Pi$
with continuous and dotted lines, respectively. The spectral bands are the boldface
vertical segments on the $\lambda$-axis.}%
\label{f3}%
\end{figure}

The sesquilinear Hermitian form on the left of (\ref{88}) is evidently
positive and closed as a consequence of Korn's inequality (\ref{9}) in the
finite cells $\omega^{0}$ and $\{  x\in\Pi:\left\vert x_{2}\right\vert
<h\}  .$ Thus, the problem (\ref{88}) is associated \cite[§ 10.1]{BiSo}
with a positive self-adjoint operator $\mathcal{A}^{\sharp}\left(
\zeta\right)  $ which has the differential expression $L^{\sharp}\left(
x,\partial_{1}+i\zeta,\partial_{2}\right)  $ and the domain%
\begin{equation}
\mathcal{D}\left(  \mathcal{A}^{\sharp}\left(  \zeta\right)  \right)
=\left\{  U^{\sharp}\in H_{\mathrm{per}\sharp}^{2}\left(  \Pi\right)
^{n}:N^{\sharp}\left(  x,\partial_{1}+i\zeta,\partial_{2}\right)  U^{\sharp
}\left(  x\right)  =0,\ \ x\in\Gamma\right\}  . \label{89}%
\end{equation}
According to Theorem \ref{thna17}, the essential spectrum of the operator
$\mathcal{A}^{\sharp}\left(  \zeta\right)  $ and therefore of problem
(\ref{88}) equals
\begin{equation}
\sigma_{es}\left(  \mathcal{A}^{\sharp}\left(  \zeta\right)  \right)
={\displaystyle\bigcup_{k\in\mathbb{N}}}B_{k}^{\sharp}\left(  \zeta\right)
,\ \ \ B_{k}^{\sharp}\left(  \zeta\right)  =\left\{  \Lambda_{k}\left(
\eta\right)  :\eta=\left(  \zeta,\eta_{2}\right)  \in\mathbf{Y}\right\}  ,
\label{90}%
\end{equation}
while according to (\ref{33})%
\begin{equation}
B_{k}^{0}={\displaystyle\bigcup_{\zeta\in\left[  0,2\pi\right)  }}%
B_{k}^{\sharp}\left(  \zeta\right)  . \label{91}%
\end{equation}
If
\begin{align}
\lambda^{\sharp}  &  \notin B_{k}^{\sharp}\left(  \zeta_{0}\right)
,\ \ \ \forall k\in\mathbb{N},\label{92}\\
\kappa^{\sharp}  &  =\dim\ker T_{0}\left(  \lambda^{\sharp};\zeta_{0}\right)
>0, \label{93}%
\end{align}
hold for some $\zeta_{0}\in\left[  0,2\pi\right)  ,$ then $\lambda^{\sharp}$
is an eigenvalue in the discrete spectrum $\sigma_{di}\left(  \mathcal{A}%
^{\sharp}\left(  \zeta_{0}\right)  \right)  $ and $\ker T_{0}\left(
\lambda^{\sharp};\zeta_{0}\right)  $ is the corresponding eigenspace. By the
continuity of functions (\ref{31}), we have%
\begin{equation}
\lambda\notin B_{k}^{\sharp}\left(  \zeta\right)  ,\ \ \ \forall
k\in\mathbb{N},\ \ \ \ \lambda\in\left(  \lambda^{\sharp}-\delta_{0}%
,\lambda^{\sharp}+\delta_{0}\right)  ,\ \ \ \ \zeta\in\left(  \zeta^{\sharp
}-\varepsilon_{0},\zeta^{\sharp}+\varepsilon_{0}\right)  , \label{94}%
\end{equation}
where $\varepsilon_{0}$ and $\delta_{0}$ are positive and the points $\zeta$
and $\zeta-2\pi$ are identified due to the evident $2\pi$-periodicity. Hence,
by general results of the perturbation theory of linear operators (cf.
\cite[Ch.\,XIII]{HiPhi}, \cite[Ch.\,9]{Kato}) the point $\left(  \zeta
^{\sharp},\lambda^{\sharp}\right)  \in\left[  0,2\pi\right)  \times
\mathbb{R}_{+}$ is the intersection of $\kappa^{\sharp}$ continuous curves
$\lambda=\lambda_{k}^{\sharp}\left(  \zeta\right)  $, which can be extended
either periodically onto the whole semi-interval $\left[  0,2\pi\right)
\ni\zeta$ or have endpoints $\left(  \zeta_{k}^{\sharp},\lambda_{k}^{\sharp
}\right)  $ at the edges of the spectral bands (\ref{90}), see Fig.\,\ref{f3}.
Thus, we have obtained an at most countable family of bounded, connected and
closed sets, that is, segments with endpoints included,

%

\begin{equation}
\mathcal{B}_{k}^{\sharp}=\left\{  \lambda=\lambda_{k}\left(  \zeta\right)
:\zeta\in\left[  \zeta_{k}^{-},\zeta_{k}^{+}\right]  \subset\left[
0,2\pi\right]  \right\}  \subset\mathbb{R}_{+}. \label{95}%
\end{equation}
The segment (\ref{95}) may be covered by a spectral band (\ref{33}%
)=(\ref{91}), or may intersect a spectral gap and in this case be contained in
the additional part (\ref{spdif}) of the essential spectrum $\sigma
_{es}\left(  \mathcal{A}\right)  $ of the problem (\ref{15}), (\ref{16}); this
fact has been mentioned in Section \ref{sect1.5} and will be of interest later.

\subsection{The regularity field for the open waveguide\label{sect3.4}}

We aim to verify the formula
\begin{equation}
\sigma_{es}\left(  \mathcal{A}\right)  =\sigma_{es}\left(  \mathcal{A}%
^{0}\right)  \cup\sigma^{\sharp}, \label{A1}%
\end{equation}
and to this end we take a point $\lambda\in\overline{\mathbb{R}_{+}}$, which
does not belong to set (\ref{A1}). Here, (\ref{A1}) contains the essential
spectrum (\ref{34}) of the purely periodic problem (\ref{3}), (\ref{4}) and
also the union (cf. formula \eqref{sigma})%
\begin{equation}
\sigma^{\sharp}=%
{\displaystyle\bigcup_{\zeta\in\left[  0,2\pi\right)  }}
\sigma_{di}\left(  \mathcal{A}^{\sharp}\left(  \zeta\right)  \right)
\label{A2}%
\end{equation}
of the discrete spectra of the operator family $\left\{  \mathcal{A}^{\sharp
}\left(  \zeta\right)  \right\}  _{\zeta\in\left[  0,2\pi\right)  },$
(\ref{89}). As explained above, the set (\ref{A2}) consists of a
union of at most countably many segments (\ref{95}).

Let the inhomogeneous problem (\ref{15}), (\ref{16}), namely%
\begin{align}
L\left(  x,\nabla\right)  u\left(  x\right)  -\lambda u\left(  x\right)   &
=f\left(  x\right)  ,\ \ \ x\in\Omega,\label{A3}\\
N\left(  x,\nabla\right)  u\left(  x\right)   &  =g\left(  x\right)
,\ \ \ x\in\partial\Omega,\nonumber
\end{align}
have the right-hand sides%
\begin{equation}
\left\{  f,g\right\}  \in L^{2}\left(  \Omega\right)  ^{n}\times
H^{1/2}\left(  \partial\Omega\right)  ^{n}. \label{A4}%
\end{equation}
We introduce a cut-off function $\mathcal{X}_{0}\in C^{\infty}\left(
\mathbb{R}^{2}\right)  $ such that%
\begin{align}
\mathcal{X}_{0}\left(  x\right)   &  =1\text{ for either }\left\vert
x_{2}\right\vert >h+d\text{, or }x_{1}<-1-d,\label{A5}\\
\mathcal{X}_{0}\left(  x\right)   &  =0\text{ for }\left\vert x_{2}\right\vert
<h\text{ and }x_{1}>-1,\ \ 0\leq\mathcal{X}_{0}\leq1.\nonumber
\end{align}
The function $\mathcal{X}_{0}$ is equal to 1 outside the semi-strip $\left(
-1-d,+\infty\right)  \times\left(  -h-d,h+d\right)  $ but vanishes inside the
smaller semi-strip $\left(  -1,+\infty\right)  \times\left(  -h,h\right)  $
which contains $\varpi^{+}$ and thus also the open waveguide. Putting%
\begin{equation}
f^{0}=\mathcal{X}_{0}f,\ \ \ \ g^{0}=\mathcal{X}_{0}g \label{A6}%
\end{equation}
gives us vector functions defined in $\Omega^{0}$ and $\partial\Omega^{0}$,
respectively. Moreover,
\begin{equation}
\left\Vert f^{0};L^{2}\left(  \Omega^{0}\right)  \right\Vert \leq\left\Vert
f;L^{2}\left(  \Omega\right)  \right\Vert ,\ \ \ \left\Vert g^{0}%
;H^{1/2}\left(  \partial\Omega^{0}\right)  \right\Vert \leq\left\Vert
g;H^{1/2}\left(  \partial\Omega\right)  \right\Vert . \label{A7}%
\end{equation}
The first inequality (\ref{A7}) is evident, while the second one is a
consequence of the following observation: in view of (\ref{A5}) and
(\ref{strip}) the function $\mathcal{X}_{0}$ equals either one or zero on each
connected component of the boundary $\partial\Omega^{0}.$
This also shows that the commutator $\left[  N,\mathcal{X}_{0}\right]  $ will
be null the in second formula of (\ref{A12}).

Since $\lambda$ stays out of the set (\ref{A1}) by assumption, the condition
(\ref{37}) is met and, according to Section \ref{sect2.3}, the problem
(\ref{37}) gets a unique solution $u^{0}\in H^{2}\left(  \Omega^{0}\right)
^{n}$ with the estimate%
\begin{align}
\left\Vert u^{0};H^{2}\left(  \Omega^{0}\right)  \right\Vert  &  \leq c\left(
\left\Vert f^{0};L^{2}\left(  \Omega^{0}\right)  \right\Vert +\left\Vert
g^{0};H^{1/2}\left(  \partial\Omega^{0}\right)  \right\Vert \right)
\leq\label{A8}\\
&  \leq c\left(  \left\Vert f;L^{2}\left(  \Omega\right)  \right\Vert
+\left\Vert g;H^{1/2}\left(  \partial\Omega\right)  \right\Vert \right)
.\nonumber
\end{align}
This solution forms the first component in the representation (notation will
be introduced later step by step)
\begin{equation}
u=\mathcal{R}\left(  \lambda\right)  \left\{  f,g\right\}  =\mathcal{X}%
_{0}u^{0}+\mathcal{X}_{\sharp}u^{\sharp}+\mathcal{X}_{b}u^{b} \label{A9}%
\end{equation}
of a parametrix for the boundary-value problem (\ref{A3}); a parametrix is by
definition a continuous operator%
\begin{equation}
\mathcal{R}\left(  \lambda\right)  :L^{2}\left(  \Omega\right)  ^{n}\times
H^{1/2}\left(  \partial\Omega\right)  ^{n}\rightarrow H^{2}\left(
\Omega\right)  ^{n} \label{A10}%
\end{equation}
such that the mapping%
\begin{equation}
\left\{  L,N\right\}  \mathcal{R}\left(  \lambda\right)  -\mathbb{I}%
:L^{2}\left(  \Omega\right)  ^{n}\times H^{1/2}\left(  \partial\Omega\right)
^{n}\rightarrow L^{2}\left(  \Omega\right)  ^{n}\times H^{1/2}\left(
\partial\Omega\right)  ^{n} \label{A11}%
\end{equation}
is compact, where $\mathbb{I}$\ stands for the identity.

We have%
\begin{align}
L\left(  \mathcal{X}_{0}u^{0}\right)   &  =\mathcal{X}_{0}f^{0}+\left[
L,\mathcal{X}_{0}\right]  u^{0}=\mathcal{X}_{0}^{2}f+\left[  L^{0}%
,\mathcal{X}_{0}\right]  u^{0},\label{A12}\\
N\left(  \mathcal{X}_{0}u^{0}\right)   &  =\mathcal{X}_{0}g^{0}=\mathcal{X}%
_{0}^{2}g ,\nonumber
\end{align}
where the term with the commutator $\left[  L^{0},\mathcal{X}_{0}\right]  $
admits the estimate%
\begin{align}
\left\Vert \left[  L,\mathcal{X}_{0}\right]  u^{0};H^{1}\left(  \Omega\right)
\right\Vert  &  =\left\Vert \left[  L^{0},\mathcal{X}_{0}\right]  u^{0}%
;H^{1}\left(  \Omega^{0}\right)  \right\Vert \leq\label{A13}\\
&  \leq c\left(  \left\Vert f;L^{2}\left(  \Omega\right)  \right\Vert
+\left\Vert g;H^{1/2}\left(  \partial\Omega\right)  \right\Vert \right)
.\nonumber
\end{align}
The difference $u^{1}=u-u^{0}$ is to be sought from the problem%
\begin{align}
L\left(  x,\nabla\right)  u^{1}\left(  x\right)   &  =f^{1}\left(  x\right)
=\left(  1-\mathcal{X}_{0}\left(  x\right)  ^{2}\right)  f\left(  x\right)
-\left[  L,\mathcal{X}_{0}\right]  u^{0}\left(  x\right)  ,\ \ \ x\in
\Omega,\label{A14}\\
N\left(  x,\nabla\right)  u^{1}\left(  x\right)   &  =g^{1}\left(  x\right)
=\left(  1-\mathcal{X}_{0}\left(  x\right)  ^{2}\right)  g\left(  x\right)
,\ \ \ x\in\partial\Omega,\nonumber
\end{align}
where the supports of the right-hand sides are located in the closed
semi-infinite strip%
\begin{equation}
\left[  -1-d,+\infty\right)  \times\left[  -h-d,h+d\right]  . \label{A15}%
\end{equation}

We now introduce the cut-off function%
\begin{equation}
\mathcal{X}_{\sharp}\left(  x\right)  =\chi\left(  x_{1}-R\right)  ,
\label{A16}%
\end{equation}
where $\chi$ and $R$ are taken from (\ref{48}) and (\ref{A}), respectively.
The products%
\begin{equation}
f^{\sharp}=\mathcal{X}_{\sharp}f^{1},\ \ g^{\sharp}=\mathcal{X}_{\sharp}g^{1}
\label{A17}%
\end{equation}
are defined in the domain $\Omega^{\sharp},$ see (\ref{53}), and its boundary
$\partial\Omega^{\sharp}.$ Moreover, the estimate
\begin{align}
\left\Vert f^{\sharp};L_{\beta}^{2}\left(  \Omega^{\sharp}\right)  \right\Vert
+\left\Vert g^{\sharp};W_{\beta}^{1/2}\left(  \partial\Omega^{\sharp}\right)
\right\Vert  &  \leq c\left(  \left\Vert f^{1};L^{2}\left(  \Omega\right)
\right\Vert +\left\Vert g^{1};H^{1/2}\left(  \partial\Omega\right)
\right\Vert \right) \label{A18}\\
&  \leq c\left(  \left\Vert f;L^{2}\left(  \Omega\right)  \right\Vert
+\left\Vert g;H^{1/2}\left(  \partial\Omega\right)  \right\Vert \right)
\nonumber
\end{align}
is valid with any weight index $\beta\in\mathbb{R}$, because $e^{\beta
\left\vert x_{2}\right\vert }\leq c_{\beta}$ in the strip (\ref{A15}), where
the supports of $f^{1}$ and $g^{1}$ are contained in.

We now make use of our assumption that $\lambda$ does not belong to the set
(\ref{A1}). As a consequence, the operator $T_{\beta}\left(  \lambda
;\zeta\right)  $ of the problem (\ref{59}), (\ref{58}) is an isomorphism
between the weighted spaces in (\ref{75}) for all
\begin{equation}
\zeta\in\left[  0,2\pi\right)  \text{ \ and \ }\beta\in\left[  0,\beta
_{0}\left(  \lambda\right)  \right)  , \label{A19}%
\end{equation}
where $\beta_{0}\left(  \lambda\right)  >0$ is determined in (\ref{82}). In
this way we apply the partial FBG-transform to the boundary-value problem
(\ref{54}) with the right-hand sides (\ref{A17}); this yields a problem of the
form (\ref{59}), (\ref{58}), and we find a unique solution $U^{\sharp}\in
W_{\beta,per\sharp}^{2}\left(  \Pi\right)  ^{n}$ for it. The above mentioned
isomorphism property of $T_{\beta}\left(  \lambda;\zeta\right)  $ guarantees
the estimate%
\begin{align}
\int_{0}^{2\pi}\left\Vert U^{\sharp};W_{\beta,per\sharp}^{2}\left(
\Pi\right)  \right\Vert ^{2}d\zeta &  \leq C\int_{0}^{2\pi}\left(  \left\Vert
F^{\sharp};L_{\beta}^{2}\left(  \Pi\right)  \right\Vert ^{2}+\left\Vert
G^{\sharp};W_{\beta}^{1/2}\left(  \Gamma\right)  \right\Vert ^{2}\right)
d\zeta\label{A20}\\
&  \leq c\left(  \left\Vert f^{\sharp};L_{\beta}^{2}\left(  \Omega^{\sharp
}\right)  \right\Vert ^{2}+\left\Vert g^{\sharp};W_{\beta}^{1/2}\left(
\partial\Omega^{\sharp}\right)  \right\Vert ^{2}\right)  ,\nonumber
\end{align}
where the latter inequality is based on (\ref{74}) (and (\ref{54})). We then
employ the inverse transform (\ref{Gin}) and obtain a solution $u^{\sharp}\in
W_{\beta}^{2}\left(  \Omega^{\sharp}\right)  ^{n}$ together with the relation%
\begin{equation}
\left\Vert u^{\sharp};H^{2}\left(  \Omega^{\sharp}\right)  \right\Vert
^{2}\leq c_{\beta}\left\Vert u^{\sharp};W_{\beta}^{2}\left(  \Omega^{\sharp
}\right)  \right\Vert ^{2}\leq C_{\beta}\int_{0}^{2\pi}\left\Vert U^{\sharp
};W_{\beta,per\sharp}^{2}\left(  \Pi\right)  \right\Vert ^{2}d\zeta.
\label{A21}%
\end{equation}
As a result, we have defined the second component in representation (\ref{A9}).

\begin{remark}
\label{remVYZ} The same cut-off function (\ref{A16}) multiplies $u^{\sharp}$
in (\ref{A9}) as well as $f^{1}$ and $g^{1}$ in (\ref{A17}). In Section
\ref{sect4.4} we explain how the second term in (\ref{A9}) must be corrected
in the case of the $\mathsf{L,V}$- and $\mathsf{Y,Z}$-shaped waveguides (see
Fig.\,\ref{f2} and Section \ref{sect1.4}).
\end{remark}

The estimates (\ref{A21}), (\ref{A22}) and (\ref{A18}) show that the product
$\mathcal{X}_{\sharp}u^{\sharp}$ in (\ref{A9}) cannot prevent the continuity
of the parametrix (\ref{A10}). Furthermore, we have%
\begin{align}
L\left(  \mathcal{X}_{\sharp}u^{\sharp}\right)   &  =\mathcal{X}_{\sharp
}f^{\sharp}+\left[  L,\mathcal{X}_{\sharp}\right]  u^{\sharp}=\mathcal{X}%
_{\sharp}^{2}f^{1}+\left[  L^{\sharp},\mathcal{X}_{\sharp}\right]  u^{\sharp
},\label{A22}\\
N\left(  \mathcal{X}_{\sharp}u^{\sharp}\right)   &  =\mathcal{X}_{\sharp
}g^{\sharp}=\mathcal{X}_{\sharp}^{2}g^{1}\nonumber
\end{align}
where%
\begin{equation}
\left\Vert \left[  L,\mathcal{X}_{\sharp}\right]  u^{\sharp};W_{\beta}%
^{1}\left(  \Omega^{\sharp}\right)  \right\Vert \leq c\left\Vert u^{\sharp
};W_{\beta}^{2}\left(  \Omega^{\sharp}\right)  \right\Vert \leq c\left(
\left\Vert f;L^{2}\left(  \Omega\right)  \right\Vert +\left\Vert
g;H^{1/2}\left(  \partial\Omega\right)  \right\Vert \right)  . \label{A23}%
\end{equation}
Since the supports of the coefficients in the
operator $\left[  L^{\sharp},\mathcal{X}_{\sharp}\right]  $ are located in the
vertical strip $\left\{  x:x_{1}\in\left[  R,R+d\right]  ,\ x_{2}\in
\mathbb{R}\right\}  ,$ the exponential weight in the norm on the left-hand
side of (\ref{A23}) helps to prove that the
mapping
\begin{equation}
L^{2}\left(  \Omega\right)  ^{n}\times H^{1/2}\left(  \partial\Omega\right)
^{n}\ni\left\{  f,g\right\}  \mapsto\left[  L^{\sharp},\mathcal{X}_{\sharp
}\right]  u^{\sharp}\in L^{2}\left(  \Omega\right)  ^{n} \label{A24}%
\end{equation}
is compact. Indeed, the mapping can be presented as a sum of a compact
operator (due to the compact embedding $H^{1}\subset L^{2}$ in a bounded
domain, the rectangle $\left[  R,R+d\right]  \times\left[  -t,t\right]  $),
and a small operator with norm $O\left(  e^{-t\beta}\right)  $ (due to the
weight which grows exponentially as $\left\vert x_{2}\right\vert >t$ and
$t\rightarrow+\infty).$

All terms with the compact embedding property, e.g. (\ref{A24}), can be
excluded from forthcoming considerations. Hence, owing to the inequality
(\ref{A13}) and formulas (\ref{A21}), (\ref{A22}), it remains to deal with the
problem (\ref{A3}) with compactly supported right-hand sides%
\[
f^{b}=\left(  1-\mathcal{X}_{0}^{2}\right)  \left(  1-\mathcal{X}_{\sharp}%
^{2}\right)  f,\ \ \ g^{b}=\left(  1-\mathcal{X}_{0}^{2}\right)  \left(
1-\mathcal{X}_{\sharp}^{2}\right)  g.
\]
Since this problem is elliptic, recall Section \ref{sect1.2}, classical
results in \cite{ADN1, ADN2} give a vector function $u^{b}\in H^{2}\left(
\Omega\right)  ^{n}$ such that%
\begin{align*}
L\left(  \mathcal{X}_{b}u^{b}\right)   &  =\mathcal{X}_{b}f^{b}+\left[
L,\mathcal{X}_{b}\right]  u^{b}=f^{b}+\left[  L,\mathcal{X}_{b}\right]
u^{b},\\
N\left(  \mathcal{X}_{b}u^{b}\right)   &  =\mathcal{X}_{b}g^{b}=g^{b},\\
\left\Vert u^{b};H^{2}\left(  \Omega\right)  \right\Vert +\left\Vert \left[
L,\mathcal{X}_{b}\right]  u^{b};H^{1}\left(  \Omega\right)  \right\Vert  &
\leq c\left(  \left\Vert f^{b};L^{2}\left(  \Omega\right)  \right\Vert
+\left\Vert g^{b};H^{1/2}\left(  \partial\Omega\right)  \right\Vert \right)
\leq\\
&  \leq c\left(  \left\Vert f;L^{2}\left(  \Omega\right)  \right\Vert
+\left\Vert g;H^{1/2}\left(  \partial\Omega\right)  \right\Vert \right),
\end{align*}
where the cut-off function $\mathcal{X}_{b}$ can be chosen as%
\[
\mathcal{X}_{b}\left(  x\right)  ={\displaystyle\prod_{\pm}}
{\displaystyle\prod_{p=1,2}} \left(  1-\chi\left(  \pm x_{1}\mp\rho\right)
\right)
\]
and $\rho\in\mathbb{N}$ is sufficiently large.

So we have constructed the parametrix $\mathcal{R}\left(  \lambda\right)  $
with all necessary properties, see (\ref{A10}) and (\ref{A11}). Since the
problem (\ref{A3}) is formally self-adjoint, we also have proved the following
assertion which, in particular, shows that the given point $\lambda$
belongs to the regularity field of the operator $\mathcal{A}$ of Section
\ref{sect1.3}.


\begin{theorem}
\label{thFRED}Assume $\lambda\in\overline{\mathbb{R}_{+}}\setminus\sigma
_{es}\left(  \mathcal{A}\right)  $, cf. (\ref{A1}). The homogeneous problem
(\ref{3}), (\ref{4}) has a finite-dimensional space $\ker\left(
\mathcal{A}-\lambda\right)  $ of solutions in $H^{2}\left(  \Omega\right)
^{n}.$ The inhomogeneous problem (\ref{A3}) admits a solution $u\in
H^{2}\left(  \Omega\right)  ^{n}$, if and only if the right-hand side
(\ref{A4}) satisfies the compatibility conditions%
\[
\left(  f,v\right)  _{\Omega}+\left(  g,v\right)  _{\partial\Omega
}=0,\ \ \ \forall v\in\ker\left(  \mathcal{A}-\lambda\right)  .
\]
This solution is defined up to an addendum in $\ker\left(  \mathcal{A}%
-\lambda\right)  $, and if in addition the orthogonality conditions%
\[
\left(  u,v\right)  _{\Omega}=0,\ \ \ \forall v\in\ker\left(  \mathcal{A}%
-\lambda\right)
\]
hold true, then it becomes unique and meets the estimate%
\[
\left\Vert u;H^{2}\left(  \Omega\right)  \right\Vert \leq c\left(  \left\Vert
f;L^{2}\left(  \Omega\right)  \right\Vert +\left\Vert g;H^{1/2}\left(
\partial\Omega\right)  \right\Vert \right)  .
\]

\end{theorem}

\subsection{The essential spectrum of the open waveguide\label{sect3.5}}

We now conclude with the main result of the paper.

\begin{theorem}
\label{thMAIN}The essential spectrum of the operator $\mathcal{A}$ of the
problem (\ref{15}), (\ref{16}) equals $\sigma_{es}\left(  \mathcal{A}\right)
=\sigma_{es}\left(  \mathcal{A}^{0}\right)  \cup\sigma^{\sharp}$, see
(\ref{A1}), where the components $\sigma_{es}\left(  \mathcal{A}^{0}\right)  $
and $\sigma^{\sharp}$ are defined in Sections \ref{sect2.2} and \ref{sect3.3}
by the formulas (\ref{34}), (\ref{A2}), respectively.
\end{theorem}

\textbf{Proof.} Thanks to Theorem \ref{thFRED}, it suffices to construct
singular Weyl sequences for the operator $\mathcal{A}$ at all points in the
set (\ref{A1}). If $\lambda\in\sigma_{es}\left(  \mathcal{A}^{0}\right)  ,$ we
may take such a sequence from Section \ref{sect2.4}, because the support
$\left[  2^{j},2^{j+1}\right]  \times\left[  2^{j},2^{j+1}\right]  $ of the
cut-off function (\ref{49}) and the entries (\ref{50}) of the sequence do not
touch the semi-strip $\left\{  x:x_{1}>0,\ \left\vert x_{2}\right\vert
<h\right\}  $, where the open waveguide lies in.

Let $\lambda$ belong to the interior of segment (\ref{95}), and recall that
the endpoints live in $\sigma_{es}\left(  \mathcal{A}^{0}\right)  $, if they
exist. By definition, there exists a $\zeta\in\left[  0,2\pi\right)  $ such
that the problem (\ref{59}), (\ref{58}) admits a non-trivial solution
$U^{\sharp}\in H_{\mathrm{per}\sharp}^{2}\left(  \Pi\right)  ^{n}$; this
generates a Floquet wave in the $x_{1}$-direction%
\begin{equation}
u^{\sharp}\left(  x\right)  =e^{i\zeta x_{1}}U^{\sharp}\left(  x\right)
\label{A31}%
\end{equation}
satisfying the homogeneous ($f^{\sharp}=0,\ g^{\sharp}=0$) problem (\ref{54})
in the periodic domain (\ref{53}). By the constructions in Sections
\ref{sect1.3} and \ref{sect3.1}, this domain coincides with $\Omega$ in the
half-plane $\left\{  x:x_{1}>0\right\}  $, where the matrices $A$ and
$A^{\sharp}$ become equal to each other, see (\ref{14}), (\ref{51}),
(\ref{opi}) and (\ref{13}), (\ref{52}). We localize the wave (\ref{A31}) by a
cut-off function similar to (\ref{49}), namely%
\begin{equation}
X_{j}^{\sharp}\left(  x\right)  =\chi_{j}\left(  x_{1}\right)  \chi\left(
x_{2}+2^{j+1}\right)  \chi\left(  2^{j+1}-x_{2}\right)  \label{A32}%
\end{equation}
where $\chi_{j}$ and $\chi$ are taken from (\ref{47}) and (\ref{48}). The
function $X_{j}^\sharp$ vanishes outside the rectangle $\left[  2^{j},2^{j+1}\right]
\times\left[  -2^{j+1},2^{j+1}\right]  $ and we will deal with indices
$j\in\mathbb{N}$ such that $2^{j}>R,$ cf. (\ref{A}). We set%
\begin{equation}
v^{\sharp j}\left(  x\right)  =X_{j}^{\sharp}\left(  x\right)  u^{\sharp
}\left(  x\right)  ,\ \ \ u^{\sharp j}\left(  x\right)  =\left\Vert v^{\sharp
j};L^{2}\left(  \Omega\right)  \right\Vert ^{-1}v^{\sharp j}\left(  x\right)
\label{A33}%
\end{equation}
and observe that, according to our choice of cut-off functions, both vector
functions in (\ref{A33}) satisfy the boundary conditions (\ref{4}).

Recall that by Theorem \ref{thna17.1} the function $u^{\sharp}\left(
x\right)  $ decays exponentially as $x_{2}\rightarrow\pm\infty$, so we can
compute the $L^{2}\left(  \Omega\right)  $-norm of $v^{\sharp j}$ as follows,
cf. (\ref{51}),
\begin{align}
\left\Vert v^{\sharp j};L^{2}\left(  \Omega\right)  \right\Vert ^{2}  &
\geq\sum_{j=2^{j}+1}^{2^{j+1}-1}\int_{0}^{1}\int_{-2^{j+1}+1}^{2^{j+1}%
-1}\left\vert U^{\sharp}\left(  x\right)  \right\vert ^{2}dx_{1}dx_{2}%
\geq\label{A34}\\
&  \geq\left(  2^{j+1}-2^{j}-2\right)  \left(  \left\Vert U^{\sharp}%
;L^{2}\left(  \Pi\right)  \right\Vert ^{2}-C\exp\left(  -2\beta2^{j+1}\right)
\right)  \geq c_{\sharp}2^{j},\ \ c_{\sharp}>0.\nonumber
\end{align}
We also obtain%
\begin{equation}
f^{\sharp j}=Lv^{\sharp j}-\lambda v^{\sharp j}=L^{\sharp}v^{\sharp j}-\lambda
v^{\sharp j}=\left[  L^{\sharp},X_{j}^{\sharp}\right]  u^{\sharp} \label{A35}%
\end{equation}
and notice that the supports of the coefficients in the commutator $\left[
L^{\sharp},X_{j}^{\sharp}\right]  $ belong to the union of two horizontal and
two vertical rectangles of sizes $d\times2^{j}$ and $2^{j+2}\times d$
respectively (they are overshadowed in Fig.\,\ref{f4},\,a). Due to the
exponential decay of $u^{\sharp},$ the horizontal rectangles only cause an
infinitesimal input into the $L^{2}\left(  \Omega\right)  $-norm of
(\ref{A35}) as $j\rightarrow+\infty,$ and the input of the vertical ones stays
uniformly bounded in $j.$ These mean that%
\[
\left\Vert \mathcal{A}u^{\sharp j}-\lambda u^{\sharp j};L^{2}\left(
\Omega\right)  \right\Vert =\left\Vert v^{\sharp j};L^{2}\left(
\Omega\right)  \right\Vert ^{-1}\left\Vert f^{\sharp j};L^{2}\left(
\Omega\right)  \right\Vert \leq c2^{-j/2}.
\]
An application of the Weyl criterion finishes the proof. $\blacksquare$

\begin{remark}
\label{remADD}As was commented in Section \ref{sect1.5}, the set (\ref{A2}) in
the representation (\ref{A1}) of $\sigma_{es}\left(  \mathcal{A}\right)  $ can
have a non-empty intersection with the essential spectrum $\sigma_{es}\left(
\mathcal{A}^{0}\right)  $ of the problem (\ref{3}), (\ref{4}) in the purely
periodic perforated plane $\Omega^{0}.$
\end{remark}

\section{Possible generalizations of the results\label{sect4}}

\subsection{Concrete problems in mathematical physics\label{sect4.1}}

$(i)$ \textit{Scalar equations.} Let $D\left(  \nabla\right)  =\nabla$ and let
$A$ be a Hermitian, positive definite $2\times2$-matrix function possessing
the properties described in Section \ref{sect1.3}. Then $n=m=2$ and%
\begin{equation}
L\left(  x,\nabla_{x}\right)  =-\nabla^{\top}A\left(  x\right)  \nabla
\label{P1}%
\end{equation}
becomes an elliptic second-order differential operator in the divergence form.
The algebraic completeness is evident, and in (\ref{8}) we have $\varrho
_{\nabla}=1$. We consider the following boundary conditions and problems: the
Neumann condition%
\begin{equation}
\nu\left(  x\right)  ^{\top}A\left(  x\right)  \nabla u\left(  x\right)
=0,\ \ \ x\in\partial\Omega, \label{P2}%
\end{equation}
which is nothing but (\ref{4}) with the co-normal derivative (\ref{7}), and
the Dirichlet boundary condition%
\begin{equation}
u\left(  x\right)  =0,\ \ \ x\in\partial\Omega, \label{P3}%
\end{equation}
for the generalized Helmholtz equation%
\begin{equation}
-\nabla^{\top}A\left(  x\right)  \nabla u\left(  x\right)  =\lambda\rho\left(
x\right)  u\left(  x\right)  ,\ \ \ x\in\Omega, \label{P4}%
\end{equation}
and the Steklov spectral problem%
\begin{align}
-\Delta u\left(  x\right)   &  =0,\ \ \ x\in\Omega,\label{P5}\\
\partial_{\nu}u\left(  x\right)   &  =\lambda u\left(  x\right)
,\ \ \ x\in\partial\Omega,\nonumber
\end{align}
where the spectral parameter $\lambda$ appears in the boundary condition while
$\partial_{\nu}=\nu\left(  x\right)  ^{\top}\nabla$ and $\Delta=\nabla^{\top
}\nabla$ denote the outward normal derivative and the Laplace operator. The
problem (\ref{P4}), (\ref{P2}) occurs in acoustics, and the problem
(\ref{P4}), (\ref{P3}) with $L=\Delta$ and $A\left(  x\right)  =\mathbb{I}$ in
the theory of quantum waveguides, see Remark \ref{remQW} and \cite{BoCa, BoCa2, CaKh, CaNaPe, CaNaRuo}. Moreover, (\ref{P5})
is related to the linear theory of water waves, cf. \cite{CaDuNa, KuMaVa} and Remark
\ref{remWW}.

$(ii)$ \textit{Elasticity.} Let $n=2,$ $m=3$ and
\begin{equation}
D\left(  \nabla\right)  =\left(
\begin{array}
[c]{ccc}%
\partial_{1} & 0 & 2^{-1/2}\partial_{2}\\
0 & \partial_{2} & 2^{-1/2}\partial_{1}%
\end{array}
\right)  ^{\top},\ \ \ \partial_{j}=\frac{\partial}{\partial x_{j}},\ \ j=1,2.
\label{P6}%
\end{equation}
Interpreting $u=\left(  u_{1},u_{2}\right)  ^{\top}$ as a displacement vector,
we employ the Voigt-Mandel matrix notation in elasticity and introduce the
strain and stress columns%
\begin{align}
{\mbox{\boldmath$\varepsilon$}}\left(  u;x\right)   &  =\left(
{\mbox{\boldmath$\varepsilon$}}_{11}\left(  u;x\right)
,{\mbox{\boldmath$\varepsilon$}}_{22}\left(  u;x\right)  ,2^{1/2}%
{\mbox{\boldmath$\varepsilon$}} _{12}\left(  u;x\right)  \right)  ^{\top
}=D\left(  \nabla\right)  u\left(  x\right)  ,\label{P7}\\
{\mbox{\boldmath$\sigma$}}\left(  u;x\right)   &  =A\left(  x\right)
{\mbox{\boldmath$\varepsilon$}}\left(  u;x\right)  =A\left(  x\right)
D\left(  \nabla\right)  u\left(  x\right)  ,\nonumber
\end{align}
which are composed from the Cartesian components
${\mbox{\boldmath$\varepsilon$}}_{jk}\left(  u\right)  =2^{-1}\left(
\partial_{j}u_{k}+\partial_{k}u_{j}\right)  $ and ${\mbox{\boldmath$\sigma$}}%
_{jk}\left(  u\right)  $ of the strain and stress tensors, respectively. Here,
$A\left(  x\right)  $ stands for the Hooke matrix of elastic moduli, and it is
real, symmetric, uniformly bounded and positive definite. The matrix $D\left(
\nabla\right)  $ in (\ref{P6}) is algebrically complete and $\varrho_{D}=2$ in
(\ref{8}), see \cite[\S \,3.7.5]{Nec} and, e.g., \cite[Example 1.12]{na262}, \cite{CaMiNa}.
We assume that $A\left(  x\right)  $ meets all requirements listed in Section
\ref{sect1.3}.

Time harmonic elastic waves with frequency $\kappa>0$ satisfy the system of
differential equations%
\begin{equation}
-\partial_{1}{\mbox{\boldmath$\sigma$}}_{j1}\left(  u;x\right)  -\partial
_{2}{\mbox{\boldmath$\sigma$}}_{j2}\left(  u;x\right)  =\kappa^{2}\rho\left(
x\right)  u_{j}\left(  x\right)  ,\ \ \ j=1,2,\ \ x\in\Omega. \label{P8}%
\end{equation}
If $\rho\left(  x\right)  =const>0$ and $\lambda=\kappa^{2},$ system
(\ref{P8}), in view of (\ref{P7}) and (\ref{P6}), takes form (\ref{3}) while
the traction-free boundary condition on $\partial\Omega$ reads as (\ref{4}).

The problem on the oscillations of a homogeneous (constant $A$), perforated
elastic plane is surely interesting for the engineering applications. To
include composite elastic materials into our considerations, we must deal with
the variable piecewise smooth Hooke matrix $A\left(  x\right)  $ and material
density $\rho\left(  x\right)  .$ Notice that the same modification can be
applied to (\ref{P4}) as well. In the next section we will show how to get rid
of the smoothness assumptions made until now when adapting our method to the
variational formulation of the elasticity problem. Finally, as a possible
application we also mention the Dirichlet problem (\ref{P3}) for the elastic
displacement vector $u$, which has in the two-dimensional case a clear
mechanical interpretation meaning that the boundaries of the holes in a thin
elastic plate are rigidly clamped.

$(iii)$ \textit{Piezoelectricity.} We set $n=3,$ $m=5$ and
\begin{equation}
D\left(  \nabla\right)  ^{\top}=\left(
\begin{array}
[c]{cc}%
D^{\mathsf{M}}\left(  \nabla\right)  ^{\top} & \mathbb{O}_{2\times2}\\
\mathbb{O}_{1\times3} & D^{\mathsf{E}}\left(  \nabla\right)  ^{\top}%
\end{array}
\right)  ,\ \ \ \mathbb{E}=\left(
\begin{array}
[c]{ccc}%
1 & 0 & 0\\
0 & 1 & 0\\
0 & 0 & 0
\end{array}
\right)  =diag\left\{  1,1,0\right\} ,  \label{P9}%
\end{equation}
where the $3\times2$-block $D^{\mathsf{M}}\left(  \nabla\right)  $ is taken
from (\ref{P6}) and $D^{\mathsf{E}}\left(  \nabla\right)  =\nabla.$ The column
$u=\left(
\begin{array}
[c]{c}%
u^{\mathsf{M}}\\
u^{\mathsf{E}}%
\end{array}
\right)  $ consists of the displacement vector $u^{\mathsf{M}}=\left(
u_{1},u_{2}\right)  ^{\top}$ and the electric potential $u^{\mathsf{E}}=u_{3}$
so that the superscripts $\mathsf{M}$ and $\mathsf{E}$ indicate mechanical and
electric fields, respectively. The so-called smart piezo-devises are able to
couple these fields of different physical nature, and this phenomenon is
described by the following system of three differential equations, see
\cite{CaNaSo, PartKud, foreign},%
\begin{equation}
D\left(  -\nabla\right)  ^{\top}A\left(  x\right)  D\left(  \nabla\right)
u\left(  x\right)  =\lambda\rho\left(  x\right)  \mathbb{E}u\left(  x\right)
,\ \ \ x\in\Omega. \label{P10}%
\end{equation}
Here, $\rho\left(  x\right)  >0$ is the material density and the matrix
$A\left(  x\right)  $ is written blockwise as%
\begin{equation}
A\left(  x\right)  =\left(
\begin{array}
[c]{cc}%
A^{\mathsf{MM}}\left(  x\right)  & A^{\mathsf{ME}}\left(  x\right) \\
A^{\mathsf{EM}}\left(  x\right)  & -A^{\mathsf{EE}}\left(  x\right)
\end{array}
\right) , \label{P11}%
\end{equation}
where $A^{\mathsf{MM}}\left(  x\right)  $ denotes the elastic Hooke matrix,
$A^{\mathsf{EE}}\left(  x\right)  $ the dielectric matrix and $A^{\mathsf{ME}%
}\left(  x\right)  =A^{\mathsf{EM}}\left(  x\right)  ^{\top}$ the
piezoelectric moduli. All matrices are real and $A^{\mathsf{MM}}\left(
x\right)  $ and $A^{\mathsf{EE}}\left(  x\right)  $ are symmetric positive
definite, hence also the matrix (\ref{P11}) is symmetric, but it is not
positive definite due to the minus sign in the bottom right-hand block. This
reflects the intrinsic transformation of the elastic energy into the electric
one and vice versa in a piezoelectric body. At the same time, the electric
potential $u^{\mathsf{E}}$ does not affect the kinetic energy at low and
middle frequencies and therefore it is absent on the right-hand side of
(\ref{P10}), cf. structure of the diagonal matrix $\mathbb{E}$ in (\ref{P9}).
We emphasize that in spite of the minus sign in (\ref{P11}) the
piezoelectricity system (\ref{P10}) is elliptic in the Douglis-Nirenberg sense
(cf. \cite[Example 1.13]{na262}).

In Section \ref{sect4.3} we demonstrate a reduction scheme from \cite{na317},
which allows us to apply the above results to the piezoelectricity system with
various boundary conditions.

$(iv)$ \textit{Plates.} Let $n=3,$ $m=6$ and
\begin{equation}
D\left(  \nabla\right)  ^{\top}=\left(
\begin{array}
[c]{cccccc}%
\partial_{1} & 0 & 2^{-1/2}\partial_{2} & 0 & 0 & 0\\
0 & \partial_{2} & 2^{-1/2}\partial_{1} & 0 & 0 & 0\\
0 & 0 & 0 & \partial_{1}^{2} & \partial_{1}^{2} & 2^{1/2}\partial_{1}%
\partial_{2}%
\end{array}
\right)  =\left(
\begin{array}
[c]{cc}%
D^{\mathsf{M}}\left(  \nabla\right)  ^{\top} & \mathbb{O}_{2\times3}\\
\mathbb{O}_{1\times3} & \nabla^{\top}D^{\mathsf{M}}\left(  \nabla\right)
^{\top}%
\end{array}
\right)  . \label{P98}%
\end{equation}
The Dirichlet problem (\ref{P10}), (\ref{P3}) with matrices (\ref{P98}) and%
\begin{equation}
A\left(  x\right)  =\left(
\begin{array}
[c]{cc}%
A^{\mathsf{MM}}\left(  x\right)  & A^{\mathsf{MB}}\left(  x\right) \\
A^{\mathsf{BM}}\left(  x\right)  & A^{\mathsf{BB}}\left(  x\right)
\end{array}
\right)  ,\ \ \ \mathbb{E}=\left(
\begin{array}
[c]{ccc}%
0 & 0 & 0\\
0 & 0 & 0\\
0 & 0 & 1
\end{array}
\right)  \label{P99}%
\end{equation}
describes the Kirchhoff model of an anisotropic inhomogeneous plate $\Omega$
with rigidly clamped boundaries of holes, which means
\begin{equation}
u \left(  x\right)  =0\in\mathbb{R}^{3},\ \ \ \partial_{\nu}u_{3}\left(
x\right)  =0,\ \ x\in\partial\Omega. \label{P97}%
\end{equation}
The vector $u=\left(  u_{1},u_{2},u_{3}\right)  ^{\top}$ includes the
longitudinal displacement $\left(  u_{1},u_{2}\right)  ^{\top}$ and the
deflection $u_{3}.$ The matrix $A\left(  x\right)  $ is real symmetric and
positive definite, but in the case of elastic symmetry it becomes
block-diagonal, i.e. $A^{\mathsf{MB}}=\left(  A^{\mathsf{BM}}\right)  ^{\top
}=\mathbb{O}_{3\times3}.$ In this case the Douglis-Nirenberg system
(\ref{P10}) decouples into two second-order equations and a fourth-order
equation. In Section \ref{sect4.3} we show that our scheme works also in the
case of high-order differential equations including the Kirchhoff plate (see also \cite{CaNaPi}).

\subsection{Operator formulation of the variational problem\label{sect4.2}}

The integral identity (\ref{17}) corresponding to the spectral problem
(\ref{15}), (\ref{16}) makes sense even in the case the matrix $A$ and scalar
$\rho$ are just bounded, measurable, and the boundary $\partial\Omega$ is
Lipschitz. It is important that all results on the solvability of model
boundary value problems in Section \ref{sect4.1} are easily adapted to their
weak formulations in the Sobolev and Kondratiev spaces (see \cite[Ch.2]{LiMa},
and for the periodic case \cite{na417}). It is straightforward to include into
our consideration any type of boundary conditions, which are covered by the
symmetric Green formula (cf. \cite[\S 2,2]{LiMa}), for example, the Dirichlet
conditions (\ref{P3}) or mixed boundary conditions.

\begin{figure}
[ptb]
\begin{center}
\includegraphics[
height=1.484in,
width=3.6772in
]%
{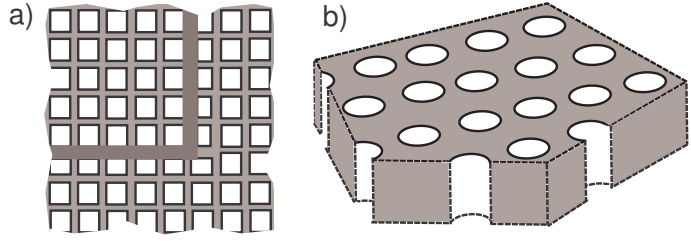}%
\caption{Other periodic geometries suitable for the present method.}%
\label{f5}%
\end{center}
\end{figure}

\begin{remark}
\label{remQW} Lipschitz domains occur for example in the grating of quantum
waveguides with long-haul thickening as in Fig.\,\ref{f5},\,a). Notice that the
domain of the operator $\mathcal{A}$ differs from $H^{2}\left(  \Omega\right)
\cap H_{0}^{1}\left(  \Omega\right)  ,$ cf. \cite{BiSk} and, e.g.,
\cite[Ch.2]{NaPl}.
\end{remark}

If the density $\rho$ is not a constant, it is useful to change the operator
formulation of the problem%
\begin{equation}
\left(  AD\left(  \nabla\right)  u,D\left(  \nabla\right)  v\right)  _{\Omega
}=\lambda\left(  \rho u,v\right)  _{\Omega},\ \ \ \forall v\in H^{1}\left(
\Omega\right)  ^{n} \label{P12}%
\end{equation}
in order to apply the theory of self-adjoint operators in Hilbert space.
Namely, having in mind the Korn inequality (\ref{AK}), we introduce in
$\mathfrak{H}=H^{1}\left(  \Omega\right)  ^{n}$ the specific scalar product
\begin{equation}
\left(  u,v\right)  _{\mathfrak{H}}=a\left(  u,v;\Omega\right)  +\left(  \rho
u,v\right)  _{\Omega}, \label{P13}%
\end{equation}
and then define the continuous, positive, symmetric, hence self-adjoint,
operator $\mathcal{S}$ in $\mathfrak{H}$ by
\begin{equation}
\left(  \mathcal{S}u,v\right)  _{\mathfrak{H}}=\left(  \rho u,v\right)
_{\Omega},\ \ \ \forall u,v\in\mathfrak{H}. \label{P14}%
\end{equation}
This turns the problem (\ref{P12}) into the abstract equation%
\begin{equation}
\mathcal{S}u=\varsigma u\text{ \ in }\mathfrak{H} \label{P15}%
\end{equation}
with the new spectral parameter%
\begin{equation}
\varsigma=\left(  1+\lambda\right)  ^{-1}. \label{P16}%
\end{equation}
The above-mentioned theory readily applies to the equation (\ref{P15}).

\begin{remark}
\label{remWW} In the perforated plane $\Omega$ (where $\vartheta^{0}$ is not
empty in (\ref{1}) and (\ref{2})) also the Steklov problem (\ref{P5}) reduces
to the equation (\ref{P15}) with parameter (\ref{P16}). According to
\cite{na449} we endow in this case the space $H^{1}\left(  \Omega\right)  $
with the scalar product $\left(  u,v\right)  _{\mathfrak{H}}=\left(  \nabla
u,\nabla v\right)  _{\Omega}+\left(  u,v\right)  _{\partial\Omega}$ and define
the "trace operator" $\mathcal{S}$ by $\left(  \mathcal{S}u,v\right)
_{\mathfrak{H}}=\left(  u,v\right)  _{\partial\Omega}.$
\end{remark}

\subsection{Reduction to integro-differential equations\label{sect4.3}}

Concerning the piezoelectricity problem with matrices (\ref{P11}), (\ref{P9})
or the plate problem with matrices (\ref{P99}), (\ref{P98}), we first mention
that the results of the papers \cite{na17, na118} can be applied here, since
they deal with general boundary value problems for Douglis-Nirenberg elliptic
systems. However, the degenerate matrix $\mathbb{E}$ on the right-hand side of
(\ref{P10}) hampers the use of the theory of self-adjoint operators in Hilbert space.

For the perforated Kirchhoff plate ($\vartheta^{0}\neq\emptyset$), the above
mentioned trick works with the new scalar product and operator $\mathcal{S}$,
\begin{align}
\left(  u,v\right)  _{\mathfrak{H}}  &  =\left(  AD\left(  \nabla\right)
u,D\left(  \nabla\right)  v\right)  _{\Omega}+\left(  \rho\mathbb{E}%
u,\mathbb{E}v\right)  _{\Omega},\label{P17}\\
\left(  \mathcal{S}u,v\right)  _{\mathfrak{H}}  &  =\left(  \rho
\mathbb{E}u,\mathbb{E}v\right)  _{\Omega},\ \ \ \forall u,v\in\mathfrak{H}%
\nonumber
\end{align}
where matrices are taken from (\ref{P99}) and (\ref{P98}). Indeed, owing to
the Dirichlet clamping condition (\ref{P97}) the bilinear form (\ref{P17}) is
a scalar product in $\mathfrak{H}=H_{0}^{1}\left(  \Omega\right)  ^{2}\times
H_{0}^{2}\left(  \partial\Omega\right)  .$

The piezoelectricity problem (\ref{P10}), (\ref{4}) requires a much more
elaborate process. Following \cite{na317}, we reduce the corresponding
variational problem%
\begin{align*}
&  \left(  A^{\mathsf{MM}}D^{\mathsf{M}}\left(  \nabla\right)  u^{\mathsf{M}%
},D^{\mathsf{M}}\left(  \nabla\right)  v^{\mathsf{M}}\right)  _{\Omega
}+\left(  A^{\mathsf{ME}}D^{\mathsf{E}}\left(  \nabla\right)  u^{\mathsf{E}%
},D^{\mathsf{M}}\left(  \nabla\right)  v^{\mathsf{M}}\right)  _{\Omega}\\
&  \ \ +\left(  A^{\mathsf{EM}}D^{\mathsf{M}}\left(  \nabla\right)
u^{\mathsf{M}},D^{\mathsf{E}}\left(  \nabla\right)  v^{\mathsf{E}}\right)
_{\Omega}-\left(  A^{\mathsf{EE}}D^{\mathsf{E}}\left(  \nabla\right)
u^{\mathsf{E}},D^{\mathsf{E}}\left(  \nabla\right)  v^{\mathsf{E}}\right)
_{\Omega}\\
&  =\left(  \rho u^{\mathsf{M}},v^{\mathsf{M}}\right)  _{\Omega},\ \ \ \forall
v\in H^{1}\left(  \Omega\right)  ^{3}%
\end{align*}
to%
\begin{equation}
\left(  A^{\mathsf{MM}}D^{\mathsf{M}}\left(  \nabla\right)  u^{\mathsf{M}%
},D^{\mathsf{M}}\left(  \nabla\right)  v^{\mathsf{M}}\right)  _{\Omega
}+E\left(  u^{\mathsf{M}},v^{\mathsf{M}}\right)  =\left(  \rho u^{\mathsf{M}%
},v^{\mathsf{M}}\right)  _{\Omega},\ \ \ \forall v^{\mathsf{M}}\in
H^{1}\left(  \Omega\right)  ^{2}, \label{P18}%
\end{equation}
where
\begin{equation}
E\left(  u^{\mathsf{M}},v^{\mathsf{M}}\right)  =\left(  A^{\mathsf{ME}%
}D^{\mathsf{E}}\left(  \nabla\right)  Ru^{\mathsf{M}},D^{\mathsf{M}}\left(
\nabla\right)  v^{\mathsf{M}}\right)  _{\Omega} \label{P19}%
\end{equation}
and $u^{\mathsf{E}}=Ru^{\mathsf{M}}\in\mathcal{H}^{1}\left(  \Omega\right)  $
is a solution of the scalar Neumann problem%
\begin{align}
&  \left(  A^{\mathsf{EE}}D^{\mathsf{E}}\left(  \nabla\right)  u^{\mathsf{E}%
},D^{\mathsf{E}}\left(  \nabla\right)  v^{\mathsf{E}}\right)  _{\Omega
}\nonumber\\
&  =\left(  A^{\mathsf{EM}}D^{\mathsf{M}}\left(  \nabla\right)  u^{\mathsf{M}%
},D^{\mathsf{E}}\left(  \nabla\right)  v^{\mathsf{E}}\right)  _{\Omega
},\ \forall v^{\mathsf{E}}\in\mathcal{H}^{1}\left(  \Omega\right)  .
\label{P20}%
\end{align}
Here, the space $\mathcal{H}^{1}\left(  \Omega\right)  $ is defined as a
completion of $C_{c}^{\infty}\left(  \overline{\Omega}\right)  $ in the norm
\[
\left(  \left\Vert \nabla u^{\mathsf{E}};L^{2}\left(  \Omega\right)
\right\Vert ^{2}+\left\Vert u^{\mathsf{E}};L^{2}\left(  K\right)  \right\Vert
^{2}\right)  ^{1/2}%
\]
and $K$ is a compactum in $\overline{\Omega}$ of positive area. Since
constants fall into $\mathcal{H}^{1}\left(  \Omega\right)  ,$ the problem
(\ref{P20}) includes one compatibility condition, which is obviously met
because the functional on the right-hand side of (\ref{P20}) degenerates on
constants. In this way, the bilinear form (\ref{P19}) is well-defined,
symmetric and positive, see \cite{na317} for details. Now the trick with a new
scalar product in $\mathfrak{H}=H^{1}\left(  \Omega\right)  ^{2}$ applies again.

Similar modifications work for other types of boundary conditions in
piezoelectricity and plates as well.

\subsection{Other geometries\label{sect4.4}}

To simplify the notation we have always used the covering of the plane with
unit squares. The cells can as well be rectangles or parallelograms, because
an affine transform does not change the crucial properties of the matrices
$D\left(  \nabla\right)  ,$ $A\left(  x\right)  $ and operators (\ref{6}),
(\ref{7}). We mention that according to \cite{na292}, the elasticity problem
in Section \ref{sect4.1} $(ii)$ preserves its matrix form under any affine
change of coordinates, though some special (non-physical!) columns of strains
and stresses will then appear. A similar procedure applies to piezoelectricity
and plates in Section \ref{sect4.1} $(iii)$ and $(iv)$.

Other types of planar coverings can be treated in a similar manner, for
instance, the diamond and honeycomb shapes, see \cite{Skri}. We do not touch
upon this generalization, which would require a total modification of the notation.

\begin{figure}
[ptb]
\begin{center}
\includegraphics[
height=1.8853in,
width=4.0127in
]%
{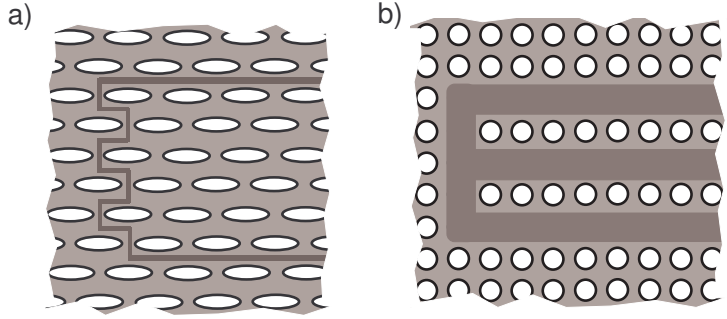}%
\caption{a) Perforation of another type. b) \textsf{E}-shaped open
waveguide.}%
\label{f6}%
\end{center}
\end{figure}

The purpose of the restriction (\ref{strip}) was solely to simplify the
definitions of cut- off functions in (\ref{47})-(\ref{50}) etc., but other
settings of holes, cf. Fig.\,\ref{f6},\,a), can be treated by our method as well.
We remark here that the homogeneous boundary conditions (\ref{4}) and
(\ref{16}) are naturally included into the definitions of the domains
(\ref{10}) and (\ref{18}) of the operators $\mathcal{A}^{0}$ and $\mathcal{A}%
$, and hence one has to take care that multiplication with the plateau
functions (\ref{49}) and (\ref{A32}) does not spoil these conditions,
otherwise the arising discrepancies must be compensated. A simplest way to
avoid the discrepancies is to construct cut-off functions, which respect the
geometry and surround all holes, like indicated in Fig.\,\ref{f6},\,a), where the
support of $\left\vert \nabla X_{j}\right\vert $ is shaded: no serious change
of calculations would follow. Of course, no modification of cut-off functions
is needed in the case of an intact periodic medium and Dirichlet conditions
(\ref{P3}).%

Our methods are sufficiently general to study the essential spectrum of the
problem (\ref{15}), (\ref{16}) in the layer $\left\{  x=\left(  y,z\right)
:y=\left(  y_{1},y_{2}\right)  \in\mathbb{R}^{2},\ \left\vert z\right\vert
<H\right\}  $ perforated periodically and perturbed inside the semi-infinite
cylinder $\mathbb{R}_{+}\times\left(  -h,h\right)  \times\left(  -H,H\right)
,$ see Fig.\,\ref{f5},\,b). As in
two-dimensional case, in this three-dimensional problem there again appear only two model spectral problems, one  in a cell and another one in an
infinite perforated prism, which can be examined with results in \cite{Skri,
Kuchbook} and \cite{na17, NaPl}. In particular our approach applies to elastic infinite
three-dimensional plates, which are perforated or have periodic bases.
On the other hand, our  approach does not yet help to examine
a problem in the space $\mathbb{R}^{3}$ with a triple-periodic perforation and with
an open periodic waveguide inside a semi-infinity circular cylinder.

\subsection{Joints of open waveguides\label{sect4.5}}

The waveguides depicted in Fig.\,\ref{f2} must keep the periodicity along all
their branches so that axes of inclined branches in Fig.\,\ref{f2},\,a),\,c),
cross the $x_{1}$-axis at angle $\phi$ with a rational $\tan\phi.$ Then each
branch gives rise to its own model problem of type (\ref{59}), (\ref{58}) in a
perforated strip $\Pi_{q}$ which is perpendicular to the branch axis and may
have width different from the main period 1; here $q=1,...,Q$ and $Q$ is the
number of branches. In this way the essential spectrum of problem (\ref{15}),
(\ref{16}) for the joint of open waveguides becomes the union of $\sigma
_{es}\left(  \mathcal{A}^{0}\right)  $ and the sets $\sigma_{1}^{\sharp
},...,\sigma_{Q}^{\sharp}$ defined in (\ref{A2}) through the discrete spectra
of the model problem in $\Pi_{1},...,\Pi_{Q}$ respectively. Theorem
\ref{thMAIN} remains true with this modification, but let us next comment the
minor changes required for its proof.

First, let us reconsider the Weyl sequence in Section \ref{sect2.4}. Since
there exists an unbounded angular opening between adjacent open waveguides, we
can choose for any $j\in N$ a point $P^{j}=\left(  P_{1}^{j},P_{2}^{j}\right)
$ such that the support of the plateau function%
\begin{equation}
X_{j}\left(  x\right)  =\chi_{j}\left(  x_{1}-P_{1}^{j}\right)  \chi
_{j}\left(  x_{2}-P_{2}^{j}\right)  \label{P77}%
\end{equation}
intersects neither the open waveguides, nor supports of $X_{1},...,X_{j-1}.$
Using (\ref{P77}), the entries (\ref{50}) of the Weyl singular will have all
the properties listed in Section \ref{sect2.4}. Similarly, in Section
\ref{sect3.5} we shift the "center" of the plateau function (\ref{A32}) along
the branch axis and so obtain the Weyl sequence elements (\ref{A33}).

Second, the construction (\ref{A9}) of the parametrix (\ref{A10}) now involves
the solutions $u_{\left(  q\right)  }^{\sharp}$ of (\ref{54}) in the domains
$\Omega_{j}^{\sharp}$ with the periodicity "cell" $\Pi_{q};$ here $q=1,...,Q.$
These solutions should be located near their branches with the help of some
cut-off functions $\mathcal{X}_{\sharp}^{q}$. However, the definition
(\ref{A16}) does not work properly, since supp$\,\mathcal{X}_{\sharp}^{q}$ may
intersect other branches of the joint and therefore new discrepancies may
appear. To avoid these discrepancies one may place the plateau function
between two neighboring branches as indicated in Fig.\,\ref{f4},\,b) by
overshading. Since supp$\,\left\vert \nabla\mathcal{X}_{\sharp}^{q}\right\vert
$ is located in a $\rho$-neighborhood of the sides of an angle domain, the
exponential weight in the Kondratiev space still leads to compact mappings of
type (\ref{A24}) and as a consequence the mapping (\ref{A11}) remains compact.
Other steps in our proof of Theorem \ref{thMAIN} remain unchanged.

We finally mention that in our notation the $\mathsf{E}$-shaped joint of open
waveguide of Fig.\,\ref{f6},\,b), only has one branch, and therefore it is directly
covered by Theorem \ref{thMAIN}.

\bigskip

\textbf{Acknowledgements}

\bigskip

G.Cardone is member of GNAMPA of INDAM. S.A. Nazarov was supported by the grant
0.38.237.2014 of St.\,Petersburg University
and the Academy of Finland project
"Mathematical approach to band-gap engineering in piezoelectric and elasticity models".
J. Taskinen was supported by the V\"ais\"al\"a Foundation of the Finnish Academy
of Sciences and Letters.

\end{document}